\documentclass[reqno]{amsart}
\usepackage{amsmath, amsthm, amssymb, amscd}

\usepackage[a4paper,width=150mm,top=25mm,bottom=25mm,centering]{geometry}

\usepackage{enumitem}
\usepackage{graphicx} 
\usepackage{xcolor}
\usepackage{comment}
\usepackage{cases}
\usepackage{array}
\usepackage{multirow}
\usepackage{longtable}
\usepackage{microtype}

\usepackage{tikz}
\usepackage{pgf}
\usepackage{tikzscale}
\usepackage{tikz-cd}
\usetikzlibrary{decorations.pathreplacing,decorations.markings}
\usetikzlibrary{arrows, knots, external, positioning}

\makeatletter
\tikzset{column sep/.code=\def\pgfmatrixcolumnsep{\pgf@matrix@xscale*(#1)},
  row sep/.code   =\def\pgfmatrixrowsep{\pgf@matrix@yscale*(#1)},
  matrix xscale/.code=%
    \pgfmathsetmacro\pgf@matrix@xscale{\pgf@matrix@xscale*(#1)},
  matrix yscale/.code=%
    \pgfmathsetmacro\pgf@matrix@yscale{\pgf@matrix@yscale*(#1)},
  matrix scale/.style={/tikz/matrix xscale={#1},/tikz/matrix yscale={#1}}}
\def\pgf@matrix@xscale{1}
\def\pgf@matrix@yscale{1}
\makeatother

\usepackage{url}
\usepackage{caption} 

\newtheorem{theorem}{Theorem}[section]

\newtheorem{lemma}[theorem]{Lemma}
\newtheorem{proposition}[theorem]{Proposition}
\newtheorem{corollary}[theorem]{Corollary}
\newtheorem{question}[theorem]{Question}

\newtheorem{claim}[theorem]{Claim}

\theoremstyle{definition}

\newtheorem{remark}[theorem]{Remark}

\renewcommand{\Bbb}{\mathbb}

\numberwithin{equation}{section}

\usepackage[bookmarksnumbered=true, colorlinks=true]{hyperref}
\AtBeginDocument{
    \hypersetup{
        citecolor=cyan,
        linkcolor=cyan,
        urlcolor=cyan
    }
}

\begin{document}

\title[Symplectic $\mathbb{Q}$-homology projective planes with rational double points]{Symplectic $\mathbb{Q}$-homology projective planes with rational double points}

\author{Woohyeok Jo}
\address{HUN-REN R\'enyi Institute of Mathematics, H-1053 Budapest, Re\'altanoda utca 13--15, Hungary}
\email{jo.woohyeok@renyi.hu}

\author{Jongil Park}
\address{Department of Mathematical Sciences and Research Institute of Mathematics, Seoul National University, Seoul 08826, Republic of Korea}
\email{jipark@snu.ac.kr}

\thanks{}
\subjclass[2020]{57K43, 57R18, 14J17, 57K41}

\keywords{Symplectic 4-orbifolds, $\mathbb{Q}$-homology projective planes, rational double points, orbifold Euler characteristic, Donaldson's diagonalization theorem}
\date\today

\begin{abstract}
We prove that a symplectic $\mathbb{Q}$-homology projective plane with only rational double points has at most five singular points and nonnegative orbifold Euler characteristic. We also show that every singularity type belongs to an explicit list of $59$ types. Of these, $58$ are realized by normal projective surfaces, while the existence of a symplectic example of type $D_4+4A_1$ is left unresolved. Thus the algebraic classification of Hwang--Keum--Ohashi extends to the symplectic category up to this single possible additional type.
\end{abstract}

\maketitle

\section{Introduction}
A normal projective complex surface whose Betti numbers agree with those of the complex projective plane $\mathbb{CP}^2$ is called an \emph{algebraic $\mathbb{Q}$-homology projective plane}, or an \emph{algebraic $\mathbb{Q}\textup{H}\mathbb{CP}^2$}. These surfaces have the smallest possible Betti numbers among normal projective complex surfaces. A nonsingular algebraic $\mathbb{Q}\textup{H}\mathbb{CP}^2$ not isomorphic to $\mathbb{CP}^2$ is called a \emph{fake projective plane}; their classification was completed through the work of Prasad--Yeung and Cartwright--Steger \cite{Prasad-Yeung-2007,Prasad-Yeung-2010,Cartwright-Steger-2010}. In the singular case, one usually restricts attention to surfaces with quotient singularities. The study of these surfaces remains active, as exemplified by the recent resolution of the algebraic Montgomery--Yang problem \cite{Jo-Park-Park-2026}. Recent works on algebraic $\mathbb{Q}\textup{H}\mathbb{CP}^2$'s include \cite{Schutt-2023,Jo-Park-Park-2024-02,Lacini-2024,Belousov-Hwang-2025,Hausen-Kiraly-Milena-2025,Jo-Park-Park-2025-01,Nagaoka-2025,Palka-Pelka-2025-02,Palka-Pelka-2025-01,Jo-Park-Park-2026-02,Liu-Liu-2026}.

A central numerical invariant in this study is the \emph{orbifold Euler characteristic}. For a normal projective complex surface $S$ with quotient singularities, it is given by
\[
e_{\mathrm{orb}}(S)=e(S)-\sum_{p\in\operatorname{Sing}(S)}\left(1-\frac{1}{|G_p|}\right),
\]
where $e(S)$ is the ordinary topological Euler characteristic and $G_p$ is the local fundamental group at $p$ \cite{Satake-1957,Kollar-2008}. Its importance is reflected in the orbifold Bogomolov--Miyaoka--Yau \textup{(BMY)} inequality, which relates the self-intersection of the canonical divisor to the topology of the surface and the orders of its local fundamental groups. We recall the following two related inequalities.

\begin{theorem}
    Let $S$ be a normal projective surface with only quotient singularities. 
    \begin{enumerate}[label=\textup{(\alph*)}]
        \item \textup{(Orbifold BMY inequality \cite{Sakai-1980,Miyaoka-1984,KoNS-1989,Megyesi-1999})}  If $K_S$ is nef, then $K^2_S \leq 3e_{\textup{orb}}(S)$.         In particular, $e_{\mathrm{orb}}(S)\geq 0$.
        \item \textup{(Weak orbifold BMY inequality \cite{KeM-1999})}  If $-K_S$ is nef, then $e_{\mathrm{orb}}(S)\geq 0$.
    \end{enumerate} 
\end{theorem}

For an algebraic $\mathbb{Q}\textup{H}\mathbb{CP}^2$ with quotient singularities, either $K_S$ or $-K_S$ is nef. The preceding theorem therefore implies $e_{\mathrm{orb}}(S)\geq0$. Since $e(S)=3$ and $|G_p|\geq2$ at every singular point, this immediately gives the preliminary bound $|\operatorname{Sing}(S)|\leq6$. Stronger restrictions follow by combining these inequalities with additional geometric and topological arguments. In particular, the following sharp general bound and its refinement for simply connected smooth loci are known.

\begin{theorem}\label{thm:algebraic_singularity_bounds}
Let $S$ be an algebraic $\mathbb Q\textup{H}\mathbb{CP}^2$ with only quotient singularities.
\begin{enumerate}[label=\textup{(\alph*)}]
\item \textup{\cite{Hwang-Keum-2011-2}} We have $|\operatorname{Sing}(S)|\leq5$, and this bound is sharp.
\item \textup{(Algebraic Montgomery--Yang problem \cite{Jo-Park-Park-2026})} If the smooth locus $S-\operatorname{Sing}(S)$ is simply-connected, then $|\operatorname{Sing}(S)|\leq3$.
\end{enumerate}
\end{theorem}

It is natural to ask to what extent these results persist in the symplectic category. We define a \emph{symplectic $\mathbb{Q}$-homology projective plane} to be a closed, effective symplectic $4$-orbifold $(X,\omega)$ with isolated singularities such that $H_*(X;\mathbb{Q})\cong H_*(\mathbb{CP}^2;\mathbb{Q})$. Every algebraic $\mathbb{Q}\textup{H}\mathbb{CP}^2$ with quotient singularities is an example after choosing an orbifold K\"ahler form.

A prominent instance in which such orbifolds play a central role is the work of Evans--Smith \cite{Evans-Smith-2018} on Lagrangian pinwheels in $\mathbb{CP}^2$. They determine which pinwheel types admit Lagrangian embeddings into $\mathbb{CP}^2$ and which collections can be embedded disjointly. These results provide a symplectic counterpart to the classification by Hacking--Prokhorov \cite{Hacking-Prokhorov-2010} of algebraic $\mathbb{Q}\textup{H}\mathbb{CP}^2$'s with quotient singularities admitting a $\mathbb{Q}$-Gorenstein smoothing to $\mathbb{CP}^2$. Related recent works on these orbifolds and symplectic constructions include \cite{Adaloglou-2025,Etnyre-2025,Chen-2026, Jo-Park-2026, Li-Ning-2026}.

A key difficulty is that no general symplectic counterpart of the orbifold BMY inequalities above is currently available. Consequently, the algebraic argument establishing the nonnegativity of the orbifold Euler characteristic cannot simply be carried over to this setting. This raises two basic questions: Does a symplectic $\mathbb{Q}\textup{H}\mathbb{CP}^2$ have nonnegative orbifold Euler characteristic, and is there a universal upper bound on its number of singular points? The two questions are closely related. As Hwang--Keum pointed out \cite[Corollary 8.4]{Hwang-Keum-2011-2}, the bound $|\operatorname{Sing}(X)|\leq 5$ extends to the symplectic setting, provided that the nonnegativity of the orbifold Euler characteristic is assumed.

In this paper, we study these questions for symplectic $\mathbb{Q}$-homology projective planes with only rational double points. This class has been extensively studied in the algebraic category \cite{Miyanishi-Zhang-1988,Miyanishi-Zhang-1993,Gurjar-Pradeep-Zhang-2002,Ye-2002,Hwang-Keum-Ohashi-2015,Schutt-2023}. Their minimal symplectic resolutions contain ADE configurations of symplectic $(-2)$-spheres, making this a natural setting in which to combine symplectic topology with lattice-theoretic obstructions. We establish the nonnegativity of the orbifold Euler characteristic and a sharp bound on the number of singular points, and then determine the possible singularity types. Our first main result is the following.

\begin{theorem}\label{thm:main_1}
    Let $X$ be a symplectic $\mathbb{Q}$-homology projective plane with only rational double points. Then $e_{\operatorname{orb}}(X)\geq 0$ and $|\operatorname{Sing}(X)|\leq 5$. Both bounds are sharp. 
\end{theorem}

The rank bound for the exceptional configurations (Proposition \ref{prop:symp_QHCP2}), together with Theorem \ref{thm:main_1}, reduces the classification to a finite list of ADE types. Combining further lattice-theoretic exclusions with known projective realizations gives the following result. 

\begin{theorem}\label{thm:classification}
    Let $X$ be a singular symplectic $\mathbb{Q}$-homology projective plane with only rational double points. Then its singularity type is one of the $59$ types listed in Table \textup{\ref{tab:classification}}. Furthermore, every type in the table, with the possible exception of $D_4+4A_1$, can be realized. In fact, each of these $58$ types is realized by a normal projective surface.
\end{theorem}

\begin{table}[!t]
\centering
\small
\renewcommand{\arraystretch}{1.4} 
\begin{tabular}{ c | >{\raggedright\arraybackslash}p{0.7\textwidth} | c }
\multicolumn{1}{c|}{Rank} & \multicolumn{1}{c|}{Singularity types} & \multicolumn{1}{c}{Number} \\
\hline \hline
1 & $A_{1}$ & 1 \\
\hline
3 & $A_{2}+A_{1}$ & 1 \\
\hline
4 & $A_{4}$ & 1 \\
\hline
5 & $D_{5}$, $A_{3}+2A_{1}$ & 2 \\
\hline
6 & $E_{6}$, $A_{5}+A_{1}$, $3A_{2}$ & 3 \\
\hline
7 & $A_{7}$, $E_{7}$, $A_{5}+A_{2}$, $D_{6}+A_{1}$, $2A_{3}+A_{1}$, $D_{4}+3A_{1}$ & 6 \\
\hline
8 & $A_{8}$, $D_{8}$, $E_{8}$, $A_{7}+A_{1}$, $2A_{4}$, $D_{5}+A_{3}$, $2D_{4}$, $E_{7}+A_{1}$, $E_{6}+A_{2}$, $A_{5}+A_{2}+A_{1}$, $D_{6}+2A_{1}$, $2A_{3}+2A_{1}$, $4A_{2}$, $\mathbf{D_{4}+4A_{1}}$ & 14 \\
\hline
9 & $A_{9}$, $D_{9}$, $A_{8}+A_{1}$, $A_{7}+A_{2}$, $A_{5}+A_{4}$, $D_{8}+A_{1}$, $D_{6}+A_{3}$, $D_{5}+A_{4}$, $D_{5}+D_{4}$, $E_{8}+A_{1}$, $E_{7}+A_{2}$, $E_{6}+A_{3}$, $A_{7}+2A_{1}$, $A_{6}+A_{2}+A_{1}$, $A_{5}+A_{3}+A_{1}$, $A_{5}+2A_{2}$, $2A_{4}+A_{1}$, $3A_{3}$, $D_{7}+2A_{1}$, $D_{6}+A_{2}+A_{1}$, $D_{5}+A_{3}+A_{1}$, $2D_{4}+A_{1}$, $E_{7}+2A_{1}$, $E_{6}+A_{2}+A_{1}$, $A_{5}+A_{2}+2A_{1}$, $A_{4}+A_{3}+2A_{1}$, $2A_{3}+A_{2}+A_{1}$, $A_{3}+3A_{2}$, $D_{6}+3A_{1}$, $D_{4}+A_{3}+2A_{1}$, $2A_{3}+3A_{1}$ & 31 
\end{tabular}
\caption{The $59$ candidate singularity types in
Theorem \ref{thm:classification}, arranged by total rank.}
\label{tab:classification}
\end{table}

The realization statement also proves the sharpness assertions in Theorem \ref{thm:main_1}: A projective example of type $2A_3+3A_1$ has five singular points and $e_{\mathrm{orb}}=0$. The $58$ projectively realized types are precisely those in the algebraic classification of Hwang--Keum--Ohashi \cite{Hwang-Keum-Ohashi-2015}, with the rank-9 realizations completed in \cite[Proposition 8.2]{Schutt-2023}. Thus, the only possible difference between the algebraic and symplectic lists is $D_4+4A_1$. See Question \ref{q:exceptional_realization} below for further discussion of this remaining type.

The classification admits sharper forms under additional homological assumptions. In what follows, $X^\circ$ denotes the smooth locus $X-\operatorname{Sing}(X)$ of $X$. We first consider a condition on the homology of the underlying space of $X$.

\begin{corollary}\label{cor:odd_homology_classification}
Let $X$ be a singular symplectic $\mathbb{Q}$-homology projective plane with only rational double points. If $H_1(X;\mathbb{Z}_2)=0$, then its singularity type is one of the $27$ types obtained by removing $D_4+4A_1$ from the entries of Table \textup{\ref{tab:classification}} of rank $<9$. Each of these $27$ types has a projective realization with a simply-connected underlying space.
\end{corollary}

Requiring the first integral homology of the smooth locus to vanish gives a substantially shorter list. See \cite[Theorem 1.7]{Jo-Park-Park-2024-02} for the algebraic counterpart. 

\begin{corollary}\label{cor:smooth_locus_homology}
Let $X$ be a singular symplectic $\mathbb{Q}$-homology projective plane with only rational double points. If $H_1(X^\circ;\mathbb{Z})=0$, then its singularity type is one of \[
A_1,\quad A_2+A_1,\quad A_4,\quad D_5,\quad
E_6,\quad E_7,\quad E_8.
\]
Each of these $7$ types has a projective realization with a simply-connected smooth locus.
\end{corollary}

\subsection*{Outline of the proofs}

We outline the proofs of Theorems \ref{thm:main_1} and \ref{thm:classification}. Let $\rho:M\to X$ be the minimal symplectic resolution, and let $R$ be the lattice generated by its exceptional spheres. If $L=\operatorname{rank}R$, then $b_2^+(M)=1$ and $K_M^2=9-L$. The adjunction formula implies $K_M$ is orthogonal to $R$. Since the real orthogonal complement of $R$ is positive definite and $1$-dimensional, this gives $L\leq9$. In the equality case $L=9$, the canonical class is torsion, and Li's results on symplectic manifolds of Kodaira dimension zero \cite{Li-2006}, together with Rohlin's theorem, identify the integral homology and the intersection lattice of $M$ with those of an Enriques surface.

The singularity bounds require more than these numerical observations. For cyclic rational double points, a key ingredient is a symplectic variant of the lattice-theoretic constraint of Jo--Park--Park \cite[Theorem 1.3]{Jo-Park-Park-2026}. This refines Donaldson's lattice-embedding obstruction by incorporating $\operatorname{spin}^c$ structures and their gluing along lens-space boundaries.

For the remaining cases, we use cyclic branched covers and symplectic blow-downs. The general pattern is \[
\begin{tikzcd}
X & M \arrow[l,"\rho"'] &
\widehat{M} \arrow[l,"\beta"'] &
Z \arrow[l,"\pi"'] \arrow[r,"f"] & Y.
\end{tikzcd}
\]
Here $\beta$ is a (possibly empty) sequence of symplectic blow-ups, $\pi$ is a cyclic cover branched over a disjoint collection of symplectic spheres, and $f$ is a sequence of symplectic blow-downs. We use Baraglia's results \cite{Baraglia-2026} to construct the cyclic branched covers, and we apply Gompf's result \cite{Gompf-1998}, which supplies symplectic forms on the branched covers together with a canonical class formula, to perform the blow-downs while tracking the canonical class. In our cases, we get 
\[
f^*K_Y=\pi^*\beta^*K_M.
\]
This identity constrains the intersection form of $Y$.

Once Theorem \ref{thm:main_1} has been established, the classification result, Theorem \ref{thm:classification}, follows the lattice-theoretic reduction used by Hwang--Keum--Ohashi \cite{Hwang-Keum-Ohashi-2015}. Their determinant and lattice-nonembedding arguments apply to the symplectic resolutions as well. The additional work is to establish the symplectic singularity bounds and to analyze the exceptional five-point configurations; the known projective examples then supply all $58$ realizations.

\smallskip

We conclude the introduction by posing further questions. We first question whether it is possible to remove the rational-double-point assumption.

\begin{question}\label{q:weak_symplectic_bmy}
Does every symplectic $\mathbb{Q}$-homology projective plane have nonnegative orbifold Euler characteristic?
\end{question}

A related question is Gompf's conjecture, stated in \cite[Problem 4.90(c)]{K3-2026}, which asks whether every closed symplectic $4$-manifold that is not a ruled surface has nonnegative second Chern number, or equivalently nonnegative Euler characteristic. Question \ref{q:weak_symplectic_bmy} asks for an analogous nonnegativity statement for the orbifold Euler characteristic in the setting of $\mathbb{Q}$-homology projective planes.

As already noted above, an affirmative answer would yield the five-point upper bound for the number of quotient singularities. The symplectic hypothesis is essential: Among oriented smooth orbifolds, the orbifold Euler characteristic can be arbitrarily negative even when all singularities are of type $A_1$. Indeed, let $Q=S^4/\langle\tau\rangle$, where \[
S^4\subset\mathbb{R}^4\times\mathbb{R},\qquad \tau(x,t)=(-x,t).
\]
The space $Q$ is topologically the suspension of $\mathbb{RP}^3$, has the rational homology of $S^4$, and its two singularities are of type $A_1$. Let $X_m:=\mathbb{CP}^2\#mQ$, where the connected sums have been taken in the smooth loci. Then $X_m$ is a simply-connected $\mathbb{Q}$-homology $\mathbb{CP}^2$, but it has $2m$ singular points and $e_{\mathrm{orb}}(X_m)=3-m$. 

An affirmative answer to Question \ref{q:weak_symplectic_bmy} would also imply the sharp upper bound on the number of pairwise disjoint symplectically embedded standard rational balls $B_{p_i,q_i}$ in $\mathbb{CP}^2$, where $p_i\geq2$. Indeed, given four such disjoint embeddings, the construction in \cite[Section 3.1]{Evans-Smith-2018} produces a symplectic $\mathbb{Q}$-homology projective plane with four Wahl singularities of orders $p_i^2$. The integers $p_i$ are pairwise coprime by \cite[Lemma 2.14]{Evans-Smith-2018}, so \[
e_{\mathrm{orb}}=-1+\sum_{i=1}^4\frac1{p_i^2}<0,
\]
contradicting an affirmative answer to the question.

There is also a stronger numerical question. For an $\omega$-compatible almost complex structure $J$, let $K_X=-c_1^{\mathrm{orb}}(TX,J)$ be the orbifold symplectic canonical class, viewed in rational cohomology. This class does not depend on the choice of $J$. Since $H^2(X;\mathbb{R})$ is one-dimensional, it is natural to use $K_X\cdot[\omega]$ to distinguish the positive, zero, and negative cases. For an algebraic $\mathbb{Q}\textup{H}\mathbb{CP}^2$ with quotient singularities, equipped with an orbifold K\"ahler form, $K_X$ is the rational cohomology class of the algebraic canonical divisor, and the canonical divisor is nef if and only if $K_X\cdot[\omega]\geq0$. This motivates the following symplectic analogue of the orbifold BMY inequality.

\begin{question}\label{q:strong_symplectic_bmy}
Let $(X,\omega)$ be a symplectic $\mathbb{Q}$-homology projective plane. If $K_X\cdot[\omega]\geq0$, does
\[
K_X^2\leq3e_{\mathrm{orb}}(X)
\]
hold?
\end{question}

On the class of symplectic $\mathbb{Q}$-homology projective planes satisfying $K_X\cdot[\omega]\geq0$, an affirmative answer would strengthen the nonnegativity assertion in Question \ref{q:weak_symplectic_bmy}. On the other hand, a realization of the exceptional type $D_4+4A_1$ would give a counterexample. Indeed, for its minimal symplectic resolution $\rho:(M,\omega_M)\to(X,\omega)$, resolving the rational double points introduces no correction terms in the canonical class, i.e., $K_M=\rho^*K_X\in H^2(M;\mathbb{Q})$. Together with Proposition \ref{prop:exceptional_intro}, this gives $K_X^2=1$ and $K_X\cdot[\omega]=K_M\cdot[\omega_M]>0$, whereas $3e_{\mathrm{orb}}(X)=3/8$.

Next, the algebraic Montgomery--Yang problem (Theorem \ref{thm:algebraic_singularity_bounds}(b)) recalled above suggests the following symplectic analogue.
\begin{question}[Symplectic Montgomery--Yang problem]\label{q:symplectic_montgomery_yang}
Let $X$ be a symplectic $\mathbb{Q}$-homology projective plane whose smooth locus $X^\circ$ is simply-connected. Does $X$ have at most three singular points?
\end{question}

For rational double points, Corollary \ref{cor:smooth_locus_homology} gives the stronger conclusion that $|\operatorname{Sing}(X)|\leq2$ under the weaker assumption $H_1(X^\circ;\mathbb Z)=0$. Thus, the remaining question concerns more general quotient singularities.

Finally, Theorem \ref{thm:classification} reduces the comparison of the algebraic and symplectic lists to a single realization question.

\begin{question}\label{q:exceptional_realization}
Does there exist a symplectic $\mathbb{Q}$-homology projective plane with singularity type $D_4+4A_1$?
\end{question}

A negative answer would show that the algebraic and symplectic lists coincide, while a positive answer would produce a singularity type with no projective realization. In Section \ref{sec:exceptional_type}, we establish further restrictions on a hypothetical symplectic realization. On the other hand, this type is realized by an oriented smooth orbifold with simply-connected underlying space and the rational homology of $\mathbb{CP}^2$; see Remark \ref{rmk:smooth_exceptional_type}. Thus, the remaining existence question genuinely concerns the symplectic category.

\smallskip

The organization of this paper is as follows. In Section \ref{sec:preliminaries}, we review symplectic resolutions, rational double points, and the lattice-theoretic constraints, and establish the rank and canonical-class restrictions on the resolution. In Section \ref{sec:main_theorem}, we prove Theorem \ref{thm:main_1}, treating the cyclic and non-cyclic cases and developing the required branched-cover arguments. In Section \ref{sec:classification}, we prove Theorem \ref{thm:classification} and Corollaries \ref{cor:odd_homology_classification} and \ref{cor:smooth_locus_homology}. In Section \ref{sec:exceptional_type}, we study the remaining type $D_4+4A_1$, deriving restrictions on symplectic realizations and constructing an oriented smooth orbifold realization. Finally, in Appendix \ref{sec:appendix}, we prove the lattice-embedding constraint, Proposition \ref{prop:lattice_constraint}, using canonical $\operatorname{spin}^c$ structures and Donaldson's diagonalization theorem.

\subsection*{Acknowledgments}
The authors thank all members of the 4-manifold topology group at Seoul National University (SNU) for their comments throughout this work. Woohyeok Jo was supported by ERC Advanced Grant KnotSurf4d. Jongil Park was supported by the National Research Foundation of Korea (NRF) grant funded by the Korean government (RS-2024-00392067). He is also affiliated with the Research Institute of Mathematics at SNU.

\section{Preliminaries}\label{sec:preliminaries}

In this section, we review the background materials and collect the preliminary results needed in the proof of the main theorems. 

\subsection{Lattices}\label{subsec:lattices}
We first fix our convention concerning lattices.

A \emph{lattice} is a pair $(L,Q)$ consisting of a free abelian group $L$ of finite rank and a nondegenerate symmetric bilinear form $Q:L\times L\to\mathbb{Q}$. For $v,w\in L$, we write $v\cdot w:=Q(v,w)$ and $v^2:=Q(v,v)$. An \emph{embedding} (resp.\ \emph{isomorphism}) of lattices is an injective group homomorphism (resp.\ group isomorphism) that preserves the bilinear forms.

A lattice $(L,Q)$ is \emph{integral} if the image of $Q$ lies in $\mathbb{Z}$, and an integral lattice $(L,Q)$ is \emph{even} if $v^2\in 2\mathbb{Z}$ for every $v\in L$. 
For an integral lattice $(L,Q)$, put $L_{\mathbb Q}:=L\otimes_{\mathbb Z}\mathbb Q$ and extend $Q$ to a $\mathbb{Q}$-bilinear form $L_{\mathbb Q} \times L_{\mathbb Q} \to \mathbb{Q}$, retaining the same notation. The \emph{dual lattice} of an integral lattice $L$ is
\[ L^*:=\{v\in L_{\mathbb Q}:Q(v,w)\in\mathbb{Z}\text{ for every }w\in L\},
\]
equipped with the restriction of $Q$. The map $L^*\to \operatorname{Hom}(L,\mathbb{Z})$, $v\mapsto Q(v,-)$, is a canonical isomorphism. Since $L$ is integral, $L\subset L^*$. The lattice is \emph{unimodular} if $L=L^*$, equivalently if $\det Q=\pm1$, where $\det Q$ denotes the determinant of a Gram matrix of $Q$ with respect to any integral basis of $L$. The quotient $L^*/L$ is called the \emph{discriminant group} of $L$. It is a finite abelian group of order $|\det Q|$. The \emph{discriminant bilinear form} is the nondegenerate pairing \[
b_L:L^*/L\times L^*/L\longrightarrow\mathbb{Q}/\mathbb{Z},
\quad b_L(v+L,w+L)=v\cdot w\pmod{\mathbb{Z}}.
\]
If $L$ is even, its \emph{discriminant quadratic form} is \[
q_L:L^*/L\longrightarrow\mathbb{Q}/2\mathbb{Z}, \quad q_L(v+L)=v^2\pmod{2\mathbb{Z}}.
\]

A sublattice $L\subset L'$ is \emph{primitive} if $L'/L$ is torsion-free, or equivalently if 
$(L\otimes\mathbb{Q})\cap L' = L$, where the intersection is taken in $L'\otimes\mathbb{Q}$.
The lattice $(L\otimes\mathbb{Q})\cap L'$ is called the \emph{primitive closure} of $L$ in $L'$. We shall use the following standard facts. 

\begin{lemma}\label{lem:discriminant} Let $(L',Q)$ be an even integral lattice, and let $L\subset L'$ be a (nondegenerate) sublattice.
    \begin{enumerate}[label=\textup{(\alph*)}]
        \item If the index $[L':L]$ is finite, then $A:=L'/L$ is isotropic with respect to the discriminant form $q_{L}$, i.e., $q_{L}|_{A}=0$, and the discriminant group $L'^*/L'$ is isomorphic to $A^\perp/A$, where \[
        A^\perp = \{x\in L^*/L:b_{L}(a,x)=0\in \mathbb{Q}/\mathbb{Z} \text{ for every $a\in A$}\}.
        \]
        \item If $(L',Q)$ is unimodular and $L$ is primitive in $L'$, then the discriminant groups $L^*/L$ and $(L^\perp)^*/L^\perp$ are isomorphic, where $L^\perp$ is the orthogonal complement of $L$ in $L'$. 
    \end{enumerate}
\end{lemma}

We will use the following lattices and conventions. \begin{itemize}
    \item For a positive integer $n$, let $\{e_1, \dots, e_n\}$ be the standard basis for $\mathbb{Z}^n$. We denote by $-\mathbb{Z}^n$ the standard diagonal negative definite lattice structure on $\mathbb{Z}^n$, where the pairing is given by $e_i\cdot e_j = -\delta_{ij}$, with $\delta_{ij}$ being the Kronecker delta.
    \item For any nonzero integer $m$, we denote by $\langle m\rangle$ the rank-one lattice generated by a vector $v$ with $v^2=m$. In particular, the standard lattice $-\mathbb{Z}^n$ is isomorphic to the direct sum $n\langle -1\rangle =\langle -1\rangle\oplus \cdots \oplus \langle-1\rangle$. 
    \item Let $X$ be a compact, oriented $4$-manifold with $H^1(\partial X;\mathbb{Z})=0$. The intersection form \[
    Q_X:H_2(X;\mathbb{Z})/\text{tor}\times H_2(X;\mathbb{Z})/\text{tor}\longrightarrow \mathbb{Z}
    \]
    gives a lattice structure on $H_2(X;\mathbb{Z})/\text{tor}$, and we simply denote this lattice by $Q_X$. The universal coefficient theorem gives an isomorphism $H^2(X;\mathbb{Z})/\text{tor}\cong \operatorname{Hom}(H_2(X;\mathbb{Z})/\text{tor},\mathbb{Z})$, and we identify this group with the dual lattice $Q_X^*$ of $Q_X$. For $4$-manifolds with boundary, cohomological squares will always refer to this pairing. 
    \item We denote by $U$ the rank-two hyperbolic lattice with Gram matrix \[
    \begin{pmatrix}
        0 & 1 \\ 1 & 0
    \end{pmatrix}.
    \]
\end{itemize}

\subsection{Minimal symplectic resolutions}
Let $(X,\omega)$ be a closed, effective symplectic $4$-orbifold with isolated singularities. Following \cite[Definition 3.2]{Cavalcanti-Fernandez-Munoz}, we define a \emph{symplectic resolution} of $(X,\omega)$ to be a closed symplectic manifold $(\tilde{X},\tilde{\omega})$ with a continuous map $\pi:\tilde{X}\to X$ satisfying the following conditions: \begin{itemize}
    \item $\pi$ restricts to a diffeomorphism $\tilde{X}-E\to X-\operatorname{Sing}(X)$, where $E:=\pi^{-1}(\operatorname{Sing}(X))$ is the \emph{exceptional set}.
    \item The exceptional set $E$ is a union of possibly intersecting smooth symplectic submanifolds of $\tilde{X}$ of codimension $\geq 2$. 
    \item $\tilde{\omega}$ and $\pi^*\omega$ agree in the complement of a small neighborhood of $E$.
\end{itemize}

In the present setting, a symplectic resolution of $(X,\omega)$ always exists. Since we shall make use of the construction, we recall the proof given in \cite{Cavalcanti-Fernandez-Munoz}.
\begin{theorem}[{\cite[Theorem 3.3]{Cavalcanti-Fernandez-Munoz}}]\label{thm:symplectic_resolution}
    Let $(X,\omega)$ be a closed, effective symplectic $4$-orbifold with isolated singularities. Then it admits a symplectic resolution.
\end{theorem}
\begin{proof}[Sketch of proof] For each $p\in\operatorname{Sing}(X)$, the equivariant Darboux theorem identifies a neighborhood of $p$ with the quotient of a neighborhood of the origin in $(\mathbb{R}^4,\omega_{\mathrm{std}})$ by a finite linear group $G_p\subset \operatorname{Sp}(4,\mathbb{R})$ acting freely away from the origin. Choose a $G_p$-invariant $\omega_{\mathrm{std}}$-compatible linear complex structure $J_p$. After a unitary linear identification with $(\mathbb{C}^2,\omega_{\mathrm{std}},i)$ and shrinking the chart, we obtain a neighborhood $U_p\cong B_p/G_p$, where $B_p\subset\mathbb{C}^2$ is a $G_p$-invariant ball and $G_p\subset \operatorname{U}(2)$. Transporting the quotient complex structure to $U_p$ makes it a K\"ahler orbifold neighborhood of $0$ in $Z_p:=\mathbb{C}^2/G_p$.

Now choose an algebraic resolution $\mu_p:\tilde{Z}_p\to Z_p$ and set $E_p:=\mu_p^{-1}(0)$, $\tilde{U}_p:=\mu_p^{-1}(U_p)$. Replacing each $U_p$ by $\tilde{U}_p$, using the identification $\tilde{U}_p-E_p\cong U_p-\{p\}$, produces a smooth manifold $\tilde{X}$ and a resolution map $\pi:\tilde{X}\to X$. For each $p$, fix a K\"ahler form $\Omega_p$ on $\tilde{Z}_p$ and the auxiliary data used in the cutoff construction. With these choices fixed, the resulting local form depends smoothly on a sufficiently small positive parameter $\epsilon_p$. This form is K\"ahler near $E_p$, agrees with the pullback of $\omega$ outside a fixed smaller neighborhood of $E_p$, and satisfies
\[
\int_C\tilde{\omega}
=\epsilon_p\int_C\Omega_p
\]
for every exceptional curve $C\subset E_p$. Performing the construction at every singular point gives the required symplectic resolution. The neighborhoods where the form is modified may be chosen arbitrarily small. \end{proof}

\begin{remark}\label{rmk:resolution_family}
In the proof of Theorem \ref{thm:symplectic_resolution}, since the set $\operatorname{Sing}(X)$ is finite, we may use a common parameter $\epsilon_p=\epsilon$ for all $p\in \operatorname{Sing}(X)$, obtaining a smooth family $\tilde{\omega}_\epsilon$, $0<\epsilon<\epsilon_1$, on $\tilde{X}$. Observe that \[
[\tilde{\omega}_\epsilon]\longleftrightarrow \left([\omega],(\epsilon\,i_p^*[\Omega_p])_p\right)
\] 
under the decomposition \[
H^2(\tilde{X};\mathbb{R})\cong H^2(X;\mathbb{R})\oplus\bigoplus_{p\in \operatorname{Sing}(X)} H^2(E_p;\mathbb{R}),
\]
where $i_p:E_p\hookrightarrow\tilde{Z}_p$ is the inclusion. Consequently, for some fixed $\eta\in H^2(\tilde X;\mathbb R)$, we have
\begin{equation}\label{eq:resolution_family}
[\tilde{\omega}_\epsilon]=\pi^*[\omega]+\epsilon\eta \longrightarrow\pi^*[\omega] \quad\text{as}\quad \epsilon\to0.
\end{equation} 
\end{remark}

\begin{remark}\label{rmk:symplectic_singularity}
From the proof of Theorem \ref{thm:symplectic_resolution}, we have the following observation. For every $p\in\operatorname{Sing}(X)$, there exist a finite subgroup $G_p\subset \operatorname{U}(2)$ and a $G_p$-invariant ball $B_p\subset\mathbb{C}^2$ such that $(U_p,\omega|_{U_p})\cong (B_p/G_p,\omega_{\mathrm{std}})$, where $G_p$ acts freely on $B_p-\{0\}$. After shrinking the neighborhoods $U_p$ if necessary, one may choose a global $\omega$-compatible orbifold almost complex structure $J$ on $X$ such that $J|_{U_p}$ is integrable and $(U_p,\omega|_{U_p},J|_{U_p}) \cong (B_p/G_p,\omega_{\mathrm{std}},i)$ as K\"ahler orbifolds. Consequently, the germ  $(X,p)$, equipped with this induced local complex structure, is analytically isomorphic to $(\mathbb{C}^2/G_p,0)$. Moreover, the resolution may be chosen so that there exist neighborhoods $\tilde{V}$ of the exceptional set $E=\pi^{-1}(\operatorname{Sing}(X))$ and $V$ of $\operatorname{Sing}(X)$, equipped with complex structures for which the symplectic forms are K\"ahler, such that $\pi|_{\tilde{V}}:\tilde{V}\to V$ is holomorphic.
\end{remark}

From now on, for every $p\in\operatorname{Sing}(X)$, we take $\mu_p:\tilde{Z}_p\to Z_p=\mathbb{C}^2/G_p$ to be the minimal algebraic resolution, so that $E_p=\mu_p^{-1}(0)$ has the weighted dual graph of the minimal resolution of the quotient singularity $0\in \mathbb{C}^2/G_p$. We shall refer to the resulting resolution $\pi:(\tilde{X},\tilde{\omega})\to (X,\omega)$ as the \emph{minimal symplectic resolution}.

Topologically, $\tilde{X}$ is obtained from $X$ by removing neighborhoods of the singular points and replacing them with neighborhoods of the corresponding exceptional divisors, or equivalently with the plumbed manifolds prescribed by their weighted dual graphs. This is also the convention used by Chen \cite[Theorem 1.5 and Remark 1.6]{Chen-2018} in his treatment of more general $4$-dimensional symplectic orbifolds. With this convention, the diffeomorphism type and the symplectic Kodaira dimension (see below) of $\tilde{X}$ are determined by $X$. 

\subsection{Symplectic \texorpdfstring{$\mathbb{Q}$}{Q}-homology projective planes}
We define a \emph{symplectic $\mathbb{Q}$-homology projective plane} to be a closed effective symplectic $4$-orbifold $(X,\omega)$ with isolated singularities such that $H_*(X;\mathbb{Q})\cong H_*(\mathbb{CP}^2;\mathbb{Q})$. As above, after choosing the induced local complex structures, its singularities are locally modeled on surface quotient singularities $(\mathbb C^2/G,0)$, where $G\subset \operatorname{U}(2)$ is finite and acts freely on $\mathbb C^2-\{0\}$. 

A singularity in the local unitary model $(\mathbb C^2/G,0)$ is a \emph{rational double point} if $G\subset \operatorname{SU}(2)$. For other characterizations, see \cite{Durfee-1979}. Rational double points are classified by the Dynkin types $A_n$ ($n\geq1$), $D_n$ ($n\geq4$), and $E_6,E_7,E_8$, according to the weighted dual graphs of their minimal resolutions; see Table \ref{tab:rational_double_points}. We use the same symbols for the associated negative definite root lattices.  

\begin{table}[t]
    \centering
    \begin{tabular}{c|c}
      Type  &  Weighted dual graph  \\ [-1em] &\\  
      \hline  
      $A_n$ & \begin{tikzpicture}[baseline=0]
\draw (-2,0) node[circle, fill, inner sep=1.2pt, black]{};
\draw (-1,0) node[circle, fill, inner sep=1.2pt, black]{};
\draw (1,0) node[circle, fill, inner sep=1.2pt, black]{};
\draw (2,0) node[circle, fill, inner sep=1.2pt, black]{};
\draw (-2,0) node[above]{$-2$};
\draw (-1,0) node[above]{$-2$};
\draw (1,0) node[above]{$-2$};
\draw (2,0) node[above]{$-2$};
\draw (0,0) node{$\cdots$};
\draw (-2,0)--(-0.5,0) (0.5,0)--(2,0) ;
\draw [thick,decorate,decoration={mirror,brace,amplitude=3pt},xshift=0pt]
	(-2,-0.2) -- (2,-0.2) node [black,midway,xshift=0pt,yshift=-10pt] 
	{$n$};
\end{tikzpicture}
      \\ [-1em] &\\  
      \hline  
       $D_n$ &  \begin{tikzpicture}[baseline=0]
\draw (0,0) node[circle, fill, inner sep=1.2pt, black]{};
\draw (1,0) node[circle, fill, inner sep=1.2pt, black]{};
\draw (2,0) node[circle, fill, inner sep=1.2pt, black]{};
\draw (1,-0.7) node[circle, fill, inner sep=1.2pt, black]{};
\draw (4,0) node[circle, fill, inner sep=1.2pt, black]{};
\draw (0,0) node[above]{$-2$};
\draw (1,0) node[above]{$-2$};
\draw (2,0) node[above]{$-2$};
\draw (3,0) node{$\cdots$};
\draw (4,0) node[above]{$-2$};
\draw (1,-0.7) node[left]{$-2$};
\draw (0,0)--(2.5,0) (3.5,0)--(4,0) (1,0)--(1,-0.7); 
\draw [thick,decorate,decoration={mirror,brace,amplitude=3pt},xshift=0pt]
	(2,-0.2) -- (4,-0.2) node [black,midway,xshift=0pt,yshift=-10pt] 
	{$n-3$};
\end{tikzpicture}  \\ [-0.8em] &\\  
      \hline
      $E_n$ & \begin{tikzpicture}[baseline=0]
\draw (-1,0) node[circle, fill, inner sep=1.2pt, black]{};
\draw (0,0) node[circle, fill, inner sep=1.2pt, black]{};
\draw (1,0) node[circle, fill, inner sep=1.2pt, black]{};
\draw (2,0) node[circle, fill, inner sep=1.2pt, black]{};
\draw (1,-0.7) node[circle, fill, inner sep=1.2pt, black]{};
\draw (4,0) node[circle, fill, inner sep=1.2pt, black]{};
\draw (-1,0) node[above]{$-2$};
\draw (0,0) node[above]{$-2$};
\draw (1,0) node[above]{$-2$};
\draw (2,0) node[above]{$-2$};
\draw (3,0) node{$\cdots$};
\draw (4,0) node[above]{$-2$};
\draw (1,-0.7) node[left]{$-2$};
\draw (-1,0)--(2.5,0) (3.5,0)--(4,0) (1,0)--(1,-0.7); 
\draw [thick,decorate,decoration={mirror,brace,amplitude=3pt},xshift=0pt]
	(2,-0.2) -- (4,-0.2) node [black,midway,xshift=0pt,yshift=-10pt] 
	{$n-4$};
\end{tikzpicture} 
    \end{tabular}
    \caption{Rational double points.}   
    \label{tab:rational_double_points}
\end{table}

For a closed symplectic $4$-manifold $(M,\omega)$, its canonical class is $K_M=-c_1(TM,J)\in H^2(M;\mathbb{Z})$, where $J$ is any $\omega$-compatible almost complex structure. Note that the class is independent of the choice of $J$. The manifold is (symplectically) \emph{minimal} if it contains no embedded symplectic $(-1)$-sphere.

\begin{proposition}\label{prop:symp_QHCP2}
    Let $X$ be a symplectic $\mathbb{Q}$-homology projective plane with only rational double points. Then the number $L$ of the exceptional curves in its minimal symplectic resolution $(\tilde{X},\tilde{\omega})$ is at most $9$, and $K_{\tilde{X}}^2=9-L\geq 0$. If $L=9$, then $\tilde{X}$ has the integral homology of an Enriques surface, $K_{\tilde{X}}$ has order two, and $Q_{\tilde{X}}\cong U\oplus E_8$.    
\end{proposition}

Recall \cite{Li-2006} that the \emph{symplectic Kodaira dimension} of a minimal symplectic $4$-manifold $(M,\omega)$ with canonical class $K$ is defined by: \[
\kappa(M,\omega)=\begin{cases}
    -\infty & \text{if $K\cdot [\omega]<0$ or $K^2<0$}, \\
    0 & \text{if $K\cdot [\omega]=0$ and $K^2=0$}, \\
    1 & \text{if $K\cdot [\omega]>0$ and $K^2=0$}, \\
    2 & \text{if $K\cdot [\omega]>0$ and $K^2>0$}.
\end{cases}
\]
For a non-minimal symplectic $4$-manifold, it is defined to be the Kodaira dimension of any of its minimal models. For a K\"ahler surface, this definition agrees with the usual holomorphic Kodaira dimension \cite[Proposition 2.7]{Li-2006}. We will use the following result. 

\begin{theorem}[{\cite[Theorem 2.4 and Proposition 6.3]{Li-2006}}]\label{thm:Li}
Let $(M,\omega)$ be a closed, minimal symplectic $4$-manifold with $\kappa(M,\omega)=0$. Then $K_M$ is either zero or of order two. In the latter case, $M$ has the integral homology of an Enriques surface.
\end{theorem}

\begin{proof}[Proof of Proposition \ref{prop:symp_QHCP2}]
The minimal symplectic resolution $\tilde{X}$ satisfies $b_1(\tilde{X})=0$ and $b_2(\tilde{X})=1+L$, so its Euler characteristic is $e(\tilde{X})=3+L$. The subgroup $R\subset H_2(\tilde{X};\mathbb{Z})/\text{tor}$ generated by the classes of the exceptional curves forms a negative definite sublattice of $Q_{\tilde{X}}$ of rank $L$. On the other hand, since $\tilde{X}$ is a closed symplectic $4$-manifold, $1-b_1(\tilde{X})+b_2^+(\tilde{X})$ must be even. It follows that $b_2^+(\tilde{X})=1$ and $b_2^-(\tilde{X})=L$, and therefore the signature of $\tilde{X}$ is $\sigma(\tilde{X})=1-L$. We conclude that $K_{\tilde{X}}^2=2e(\tilde{X})+3\sigma(\tilde{X})=9-L$. 

Next, the adjunction formula implies that $K_{\tilde{X}}$ is orthogonal to the exceptional curves, so $K_{\tilde{X}}$ is contained in the orthogonal complement $R^\perp$ of $R$, which is a 1-dimensional positive definite subspace. This yields $K_{\tilde{X}}^2\geq 0$, and hence $L\leq 9$.

To prove the final statement, suppose $L=9$ so that $K_{\tilde{X}}^2=0$. Since $K_{\tilde{X}}$ is contained in the positive definite subspace $R^\perp$, we see that $K_{\tilde{X}}$ is torsion in $H^2(\tilde{X};\mathbb{Z})$. In particular, $K_{\tilde{X}}\cdot [\tilde{\omega}]=0$ and $K_{\tilde{X}}^2=0$. Moreover, there cannot exist a symplectic $(-1)$-sphere in $\tilde{X}$ by the adjunction formula; hence $(\tilde{X},\tilde{\omega})$ is minimal. Therefore, $\kappa(\tilde{X},\tilde{\omega})=0$ and Theorem \ref{thm:Li} implies that $K_{\tilde{X}}$ is either trivial or of order two. But $K_{\tilde{X}}$ must be nontrivial, since otherwise $\tilde{X}$ would be spin with $\sigma(\tilde{X})=-8$, contradicting Rohlin's theorem \cite{Rohlin-1952}. Hence $K_{\tilde{X}}$ must be of order two, and $\tilde{X}$ has the integral homology of an Enriques surface. The fact $2K_{\tilde{X}}=0$ implies that $Q_{\tilde{X}}$ is an even lattice. Since $b_2^+(\tilde{X})=1$ and $b_2^-(\tilde{X})=9$, it follows that $Q_{\tilde{X}}\cong U\oplus E_8$. 
\end{proof}

\begin{proposition}\label{prop:determinant_square}
    Let $X$ be a symplectic $\mathbb{Q}$-homology projective plane with only rational double points. Suppose $L<9$, where $L$ is the number of the exceptional curves in the minimal symplectic resolution $\tilde{X}$. Then $(9-L)|\det R|$ is a square of a positive integer, where $R\subset Q_{\tilde{X}}$ is the sublattice generated by the classes of the exceptional curves. Moreover, $Q_{\tilde{X}}\cong \langle +1\rangle \oplus L\langle -1\rangle$ if $L\geq 2$. 
\end{proposition}

\begin{proof}
    Let $k\in H_2(\tilde{X};\mathbb{Z})/\text{tor}$ denote the Poincar\'e dual of $K_{\tilde{X}}$. We have seen in the proof of Proposition \ref{prop:symp_QHCP2} that $k^2=9-L>0$ and that $k\in R^\perp$. Therefore, $R\oplus \mathbb{Z}\cdot k$ is a full-rank sublattice of the unimodular lattice $Q_{\tilde{X}}$, and therefore $|\det (R\oplus \mathbb{Z}\cdot k)|=|\det R|\cdot (9-L)$ is the square of the index $[Q_{\tilde{X}}:R\oplus \mathbb{Z}\cdot k]$. The last assertion follows from the fact that $Q_{\tilde{X}}$ is an odd unimodular lattice with signature $(1,L)$. 
\end{proof}

\section{Proof of Theorem \ref{thm:main_1}}\label{sec:main_theorem}

In this section, we give a proof of Theorem \ref{thm:main_1}. We divide the proof into two cases according to whether all rational double points are of type $A$, or at least one of them is of type $D$ or $E$. In both cases, one of our main tools is the following variant of the lattice-theoretic constraint of \cite[Theorem 1.3]{Jo-Park-Park-2026}, whose proof is given in Appendix \ref{sec:appendix}.

\begin{proposition}\label{prop:lattice_constraint}
    Let $M$ be a closed, connected, symplectic $4$-manifold with $b_1(M)=0$ and $b_2^+(M)=1$. Suppose $M$ contains pairwise disjoint symplectic configurations of types $A_{k_1},\dots,A_{k_n}$, where $k_1+\cdots+k_n=b_2^-(M)$. Then there exists a codimension-one lattice embedding \[
    \iota:\bigoplus_{i=1}^n \langle -(k_i+1)\rangle \longrightarrow -\mathbb{Z}^{n+1},
    \]
    such that, for a chosen generator $v_i$ of the $i$-th summand $\langle-(k_i+1)\rangle$, \[
    \mathbf{K}\cdot \iota(v_i)=k_i+1 \quad (i=1,\dots,n),
    \]
    where $\mathbf{K}=3e_1+e_2+\cdots+e_{n+1}\in -\mathbb{Z}^{n+1}$.
\end{proposition}

\subsection{The cyclic case}

In this subsection, we prove Theorem \ref{thm:main_1} for symplectic $\mathbb{Q}$-homology projective planes with only cyclic rational double points. In this case, the minimal symplectic resolution satisfies the assumptions of Proposition \ref{prop:lattice_constraint} above.

For a vector $w\in -\mathbb{Z}^n$, let us define its \emph{support} to be the set $\{i\in \{1,\dots,n\}: e_i\cdot w\neq 0\}$. Thus, if we write $w=\sum_{i=1}^n w_ie_i$, then its support consists of those $i$'s with $w_i\neq 0$. 

\begin{proposition}\label{prop:cyclic_case}
     Let $X$ be a symplectic $\mathbb{Q}$-homology projective plane with only cyclic rational double points. If $k$ is the number of singularities of type $A_1$ and $n=|\operatorname{Sing}(X)|$, then \[
    k\leq \left\lfloor \frac{n+1}{2} \right\rfloor
    \]
    and $n\leq 5$.
\end{proposition}

\begin{proof} 
Apply Proposition \ref{prop:lattice_constraint} to the minimal symplectic resolution of $X$. Observe that if a vector $u\in-\mathbb{Z}^{n+1}$ satisfies $u^2=-2$ and $\mathbf{K}\cdot u=2$, then it must be of the form $u=-e_1+e_a$ for some $a>1$ or $-e_a-e_b$ for some $a,b>1$ with $a\neq b$. Since there are $k$ singular points of type $A_1$, there must exist pairwise orthogonal vectors $u_1,\dots,u_k\in -\mathbb{Z}^{n+1}$ satisfying $u^2_i=-2$ and $\mathbf{K}\cdot u_i=2$. There can exist at most one $u_i$ of the form $-e_1+e_a$, and the supports of the $u_i$'s must be pairwise disjoint. It follows that $k\leq \lfloor \frac{n+1}{2} \rfloor$. 

Now, let $L$ be the number of the exceptional curves in the minimal symplectic resolution. Then $L\geq k+2(n-k)=2n-k$. Since $L\leq 9$ by Proposition \ref{prop:symp_QHCP2}, we obtain \[
2n-9 \leq k\leq \left\lfloor \frac{n+1}{2} \right\rfloor,
\]
which implies $n\leq 6$. To show that $n\leq 5$, suppose $n=6$. Then $X$ must have exactly $3$ singular points of type $A_1$ and $3$ of type $A_2$. By Proposition \ref{prop:lattice_constraint}, there exist pairwise orthogonal vectors $u_1,u_2,u_3,v_1,v_2,v_3\in -\mathbb{Z}^7$ such that $u_i^2=-2$, $v_i^2=-3$, $\mathbf{K}\cdot u_i=2$, and $\mathbf{K}\cdot v_i=3$. In particular, each $v_i$ must be of the form $-e_1+e_a-e_b$ for some $a,b>1$ with $a\neq b$, or $-e_a-e_b-e_c$ for some distinct $a,b,c>1$.

First, consider the case where each $u_i$ is of the form $-e_a-e_b$ with $a,b>1$ and $a\neq b$. Then by permuting the vectors $e_2,\dots,e_7$ (note that this does not change $\mathbf{K}$) if necessary, we may assume $u_1=-e_2-e_3$, $u_2=-e_4-e_5$, and $u_3=-e_6-e_7$. Since $v_i\cdot u_j=0$ for all $i,j$, it follows that each $v_i$ must be of the form $-e_1+e_a-e_b$, where $\{a,b\}$ is one of $\{2,3\}$, $\{4,5\}$, or $\{6,7\}$. But then we cannot have $v_i\cdot v_j=0$ for $i\neq j$, which is a contradiction. Next, consider the case where some $u_i$ is of the form $-e_1+e_a$. In this case, we may assume without loss of generality that $u_1=-e_1+e_2$, $u_2=-e_3-e_4$, and $u_3=-e_5-e_6$. But then each $v_i$ must be $-e_1+e_7-e_2$, so again we cannot have $v_i\cdot v_j=0$ for $i\neq j$. We conclude that the case $n=6$ cannot occur. 
\end{proof} 

The next goal is to show that if $X$ is a symplectic $\mathbb{Q}$-homology projective plane with only cyclic rational double points, then $e_{\text{orb}}(X)\geq 0$ (Corollary \ref{cor:e_orb_cyclic}). To this end, we need to eliminate the cases $A_3+2A_2+2A_1$ (this notation means that $X$ has one singular point of type $A_3$, two of type $A_2$, and two of type $A_1$) and $4A_2+A_1$. 

\begin{proposition}\label{prop:A_32A_22A_1}
There is no symplectic $\mathbb{Q}$-homology projective plane with singularity type $A_3+2A_2+2A_1$. 
\end{proposition}

\begin{proof}
    We argue as in the proof of Proposition \ref{prop:cyclic_case}. Suppose there exists a symplectic $\mathbb{Q}$-homology projective plane with singularity type $A_3+2A_2+2A_1$. By applying Proposition \ref{prop:lattice_constraint} to the minimal symplectic resolution, we see that there exist pairwise orthogonal vectors $u_1,u_2,v_1,v_2,w\in -\mathbb{Z}^6$ such that $u_i^2=-2$, $v_i^2=-3$, $w^2=-4$, $\mathbf{K}\cdot u_i=2$, $\mathbf{K}\cdot v_i=3$, and $\mathbf{K}\cdot w=4$. We will not make use of the vector $w$. 

    Recall from the proof of Proposition \ref{prop:cyclic_case} that each $u_i$ is either of the form $-e_1+e_a$ ($a>1$) or $-e_a-e_b$ ($a,b>1$, $a\neq b$), and that each $v_i$ is either of the form $-e_1+e_a-e_b$ ($a,b>1$, $a\neq b$) or $-e_a-e_b-e_c$ ($a,b,c>1$ distinct). Thus, without loss of generality, we may assume either $(u_1,u_2)=(-e_2-e_3, -e_4-e_5)$ or $(u_1,u_2)=(-e_1+e_2, -e_3-e_4)$. In the former case, since $v_i\cdot u_j=0$, both $v_i$ must be of the form $-e_1+e_a-e_b$ where $\{a,b\}$ is either $\{2,3\}$ or $\{4,5\}$, but then we cannot have $v_1\cdot v_2=0$. In the latter case, both $v_i$ must be either $-e_1+e_5-e_2$ or $-e_1+e_6-e_2$, so again we cannot have $v_1\cdot v_2=0$. This completes the proof. 
\end{proof}

Unlike $A_3+2A_2+2A_1$, the type $4A_2+A_1$ satisfies the lattice constraint of Proposition \ref{prop:lattice_constraint}. Thus, this case requires a rather different geometric argument. 

\begin{proposition}\label{prop:4A_2A_1}
    There is no symplectic $\mathbb{Q}$-homology projective plane with singularity type $4A_2+A_1$.
\end{proposition}

\begin{corollary}\label{cor:e_orb_cyclic}
    If $X$ is a symplectic $\mathbb{Q}$-homology projective plane with only cyclic rational double points, then $e_{\operatorname{orb}}(X)\geq 0$.
\end{corollary}

\begin{proof}
    Let $n$ be the number of singular points, and suppose that the singular points are of types $A_{k_1},\dots,A_{k_n}$ with $k_1\leq \cdots \leq k_n$. Then the orbifold Euler characteristic of $X$ is given by \[
    e_{\textup{orb}}(X)=3-n+\sum_{i=1}^n \frac{1}{k_i+1}.
    \]
    Also, we have $n\leq 5$ by Proposition \ref{prop:cyclic_case}, and $\sum_{i=1}^n k_i\leq 9$ by Proposition \ref{prop:symp_QHCP2}. If $n\leq 3$, then clearly $e_{\operatorname{orb}}(X)>0$. If $n=4$, then $\sum_i k_i\leq 9$ implies 
    $\sum_i \frac{1}{(k_i+1)} \geq \frac{16}{\sum_i (k_i+1)} \geq \frac{16}{13}$ 
    by Cauchy--Schwarz, so again $e_{\operatorname{orb}}(X)> 0$. 
    Now suppose $n=5$ and $e_{\operatorname{orb}}(X)<0$. Then 
    \[
    \sum_{i=1}^5 \frac{1}{k_i+1}<2 \quad \text{and} \quad \sum_{i=1}^5 (k_i+1)\leq 14,
    \]
    so we must have either $(k_1,\dots,k_5)=(1,1,2,2,3)$ or $(k_1,\dots,k_5)=(1,2,2,2,2)$. But both cases are impossible by Propositions \ref{prop:A_32A_22A_1} and \ref{prop:4A_2A_1}, respectively. Hence we must have $e_{\operatorname{orb}}(X)\geq 0$ even if $n=5$.     
\end{proof}

It remains to prove Proposition \ref{prop:4A_2A_1}. We will use the following lemmas.

\begin{lemma}[{\cite[Lemma 4.1 and Proposition 4.4]{Baraglia-2026}}]\label{lem:Baraglia}
    Let $M$ be a closed, connected, oriented, smooth $4$-manifold, and let $S_1,\dots, S_r$, $r\geq 1$, be pairwise disjoint, oriented, smoothly embedded spheres of negative self-intersection. For a prime $p$, consider the following conditions. \begin{enumerate}[label=\textup{(\alph*)}]
        \item $H_1(M;\mathbb{Z})$ is finite with order coprime to $p$.
        \item Any $r-1$ of $[S_1],\dots,[S_r]$ are linearly independent in $H_2(M;\mathbb{Z}_p)$.
        \item $\sum_{i=1}^r a_i [S_i]=0\in H_2(M;\mathbb{Z}_p)$ for some $a_1,\dots,a_r\in \mathbb{Z}_p^*$.
    \end{enumerate}
    If \textup{(c)} holds, then there exists a connected cyclic branched cover $\tilde{M}\to M$ of degree $p$ with branch locus $\bigcup_{i=1}^r S_i$, whose monodromy sends a positively oriented meridian of $S_i$ to $a_i\in\mathbb{Z}_p$. If all three conditions hold, then $b_1(\tilde{M})=0$.  
\end{lemma}

In this section, we will use only the existence statement of Lemma \ref{lem:Baraglia}. The assertion $b_1(\tilde{M})=0$ will be used in Section \ref{sec:exceptional_type}.

\begin{lemma}[{\cite[Lemma 1]{Gompf-1998}}]\label{lem:Gompf}
Let $(M,\omega)$ be a closed symplectic $4$-manifold, and let $B\subset M$ be a closed, possibly disconnected, embedded symplectic surface. Suppose $\pi:\tilde{M}\to M$ is a smooth branched covering with branch locus $B$, and $\pi^{-1}(B)$ has components $B_1,\dots,B_s$, oriented by pulling back the symplectic orientation of $B$. Let $d_i\geq 1$ denote the multiplicity of $\pi$ along $B_i$. Then $\tilde{M}$ admits a symplectic form $\tilde{\omega}$ such that $[\tilde{\omega}]=\pi^*[\omega]\in H^2(\tilde{M};\mathbb{R})$, \[
c_1(\tilde{M},\tilde{\omega})= \pi^*c_1(M,\omega) + \sum_{i=1}^s (1-d_i)\operatorname{PD}[B_i]\in H^2(\tilde{M};\mathbb{Z}),
\]
and each $B_i$ is a $\tilde{\omega}$-symplectic submanifold with $\tilde{\omega}|_{B_i}=(\pi|_{B_i})^*(\omega|_{\pi(B_i)})$. 
\end{lemma}

Since we shall slightly modify the construction, we recall the main steps of the proof:

\begin{proof}[Sketch of proof]
The pullback form $\pi^*\omega$ is symplectic away from the ramification locus, but it degenerates along every component $B_i$ with $d_i\geq2$. Choose a closed symplectic tubular neighborhood $N$ of $B$ equipped with the Hamiltonian circle action which rotates the normal fibers, and let $H:N\to [0,\epsilon]$ be its Hamiltonian, normalized so that $B=H^{-1}(0)$. The punctured neighborhood $N-B$ can be written as a
Hamiltonian collar, and symplectic cutting at the level $0$ recovers the original symplectic neighborhood $N$. For each component $B_i\subset\pi^{-1}(B)$, pull this collar and its circle action back to a punctured tubular neighborhood $N_i-B_i$. The lifted circle action has Hamiltonian $d_iH\circ\pi$. Symplectic cutting at the level $0$ then produces a symplectic form $\tilde{\omega}$ on $\tilde{M}$ which agrees with
$\pi^*\omega$ outside the neighborhoods $N_i$ and satisfies $\tilde{\omega}|_{B_i}=(\pi|_{B_i})^*(\omega|_{\pi(B_i)})$. 

To determine its cohomology class, represent an arbitrary $2$-dimensional homology class by a surface which meets each $N_i$ in a union of invariant normal disks. The construction gives $\int_D\tilde{\omega} = d_i\int_{\pi(D)}\omega$ on each such disk $D$, while $\tilde{\omega}=\pi^*\omega$ outside the $N_i$. It follows that $[\tilde{\omega}]=\pi^*[\omega]$.

For the computation of $c_1(\tilde{M},\tilde{\omega})$, one chooses an $\omega$-tame almost complex structure $J$ on $M$ such that $B$ and the normal disk fibers of $N$ are $J$-holomorphic. As explained in \cite{Gompf-1998}, $J$ lifts to a continuous almost complex structure $\tilde{J}$ on $\tilde{M}$ tamed by $\tilde{\omega}$. The remaining computation of $c_1(\tilde{M},\tilde{\omega})$ uses this $\tilde{J}$; we omit the details.
\end{proof}

\begin{proof}[Proof of Proposition \ref{prop:4A_2A_1}] 
The proof can be regarded as a symplectic analog of the proof of the corresponding result \cite[Proposition 1.4]{Schutt-2019} in algebraic geometry, which states that there is no Enriques surface containing a $4A_2+A_1$-configuration of rational curves.

    Suppose $X$ is a symplectic $\mathbb{Q}$-homology projective plane with singularity type $4A_2+A_1$, and let $(M,\omega_M)\to X$ denote the minimal symplectic resolution. By Proposition \ref{prop:symp_QHCP2}, $M$ has the integral homology of an Enriques surface, $K_M\in H^2(M;\mathbb{Z})$ is of order two, and $Q_M\cong U\oplus E_8$. In particular, $H_1(M;\mathbb{Z})\cong \mathbb{Z}_2$ and $\operatorname{tor}H_2(M;\mathbb{Z})\cong \mathbb{Z}_2$. 
    
    Let $C_i\cup C_i'$ ($i=1,2,3,4$) denote the four $A_2$-configurations on $M$, so that $C_i, C_i'$ are embedded symplectic $(-2)$-spheres with $C_i\cdot C_i'=1$. Also, let $A$ denote the $A_1$-sphere in $M$. Let $R\subset L_M:= H_2(M;\mathbb{Z})/\text{tor}$ denote the negative definite sublattice of $Q_M$ generated by the classes of $C_i,C_i',A$. 

      \begin{claim}\label{claim1}
      After relabeling the four $A_2$-configurations and interchanging $C_i$ and $C_i'$ in some of them, the class $D=\sum_{i=1}^3 ([C_i]-[C_i'])$ is $3$-divisible in $H_2(M;\mathbb{Z})$.
      \end{claim}
      
\begin{proof}[Proof of Claim \ref{claim1}]
Since multiplication by $3$ is an automorphism of $\operatorname{tor}H_2(M;\mathbb{Z})\cong \mathbb{Z}_2$, it suffices to prove $3$-divisibility in $L_M=H_2(M;\mathbb{Z})/\text{tor}$. For $i=1,\dots,4$, put $\gamma_i:=([C_i]-[C_i'])/3\in R\otimes \mathbb{Q}$. Then $\gamma_i\cdot [C_i]=-1$, $\gamma_i\cdot [C_i']=1$, and $\gamma_i^2=-2/3$. In particular, $\gamma_i\in R^*$, and the classes $\gamma_i+R$ form a $\mathbb{Z}_3$-basis for the $3$-primary part $(R^*/R)_3\cong \mathbb{Z}_3^{\oplus 4}$ of the discriminant group $R^*/R\cong \mathbb{Z}_3^{\oplus 4}\oplus \mathbb{Z}_2$. 

Let $\bar{R}:=(R\otimes \mathbb{Q})\cap L_M$ denote the primitive closure of $R$ in $L_M$. By Lemma \ref{lem:discriminant}(a), $\bar{R}/R$ is isotropic with respect to $q_R$. Also, since $\bar{R}$ is primitive of codimension one in the unimodular lattice $L_M$, its discriminant group is cyclic by Lemma \ref{lem:discriminant}(b). On the other hand, since $3^4\cdot 2=|\det R|=[\bar{R}:R]^2\cdot |\det \bar{R}|$, we have $[\bar{R}:R]\in \{1,3,9\}$. But the index cannot be $1$ because $R^*/R$ is not cyclic. Hence $\bar{R}/R$ is a nonzero subgroup of $(R^*/R)_3$. Choose a nonzero element  
\[
h=\sum_{i=1}^4a_i\gamma_i+R\in \bar{R}/R, \quad a_i\in\{0,\pm1\}.
\]
If $w$ is the number of nonzero coefficients, then $q_R(h)=-2w/3 \pmod{2\mathbb{Z}}$. Since $\bar{R}/R$ is isotropic with respect to $q_R$, we must have $w=3$. Thus, after relabeling the configurations and interchanging some $C_i$ and $C_i'$, we may assume that $\gamma_1+\gamma_2+\gamma_3\in \bar{R}\subset L_M$. Since $3(\gamma_1+\gamma_2+\gamma_3)=D$, we conclude that $D$ is $3$-divisible in $L_M$. 
\end{proof}

    Let $\beta:(\widehat{M},\omega_{\widehat{M}})\to M$ denote the symplectic blow-up at the intersection points $C_i\cap C_i'$ for $i=1,2,3$. Let $\widehat{C}_i$ and $\widehat{C}_i'$ denote the proper transforms of $C_i$ and $C_i'$, respectively. These are pairwise disjoint symplectic $(-3)$-spheres in $\widehat{M}$. Let $E_i$ denote the exceptional sphere over $C_i\cap C_i'$. Identifying $H_2(\widehat{M};\mathbb{Z})$ with $H_2(M;\mathbb{Z})\oplus \bigoplus_{i=1}^3 \mathbb{Z}[E_i]$, we have $[\widehat{C}_i]=[C_i]-[E_i]$ and $[\widehat{C}_i']=[C_i']-[E_i]$, so it follows that $\widehat{D}:=\sum_{i=1}^3 ([\widehat{C}_i]-[\widehat{C}_i'])$ is $3$-divisible in $H_2(\widehat{M};\mathbb{Z})$. Hence \[
    \sum_{i=1}^3 \left([\widehat{C}_i]+2[\widehat{C}_i'] \right)=0\in H_2(\widehat{M};\mathbb{Z}_3).
    \]

  By Lemma \ref{lem:Baraglia}, there is a connected cyclic branched cover $\pi:Z\to \widehat{M}$ of degree $3$ with branch locus $\bigcup_{i=1}^3(\widehat{C}_i\cup \widehat{C}_i')$, whose monodromy sends positively oriented meridians of $\widehat C_i$ and $\widehat C_i'$ to $1$ and $2$ in $\mathbb{Z}_3$, respectively. By Lemma \ref{lem:Gompf}, we can choose a symplectic form $\omega_Z$ on $Z$ with $[\omega_Z]=\pi^*[\omega_{\widehat{M}}]$ and $K_Z=\pi^*K_{\widehat{M}}+2\sum_{i=1}^3(\operatorname{PD}[T_i]+\operatorname{PD}[T_i'])$, where $T_i=\pi^{-1}(\widehat{C}_i)$ and $T_i'=\pi^{-1}(\widehat{C}_i')$ are symplectic surfaces in $Z$. We can further assume that $F_i:=\pi^{-1}(E_i)$ ($i=1,2,3$) are also symplectic surfaces in $Z$. Indeed, recall that the original $A_2$-configurations in $M$ lie in K\"ahler neighborhoods (Remark \ref{rmk:symplectic_singularity}). Perform the three blow-ups holomorphically in the local K\"ahler neighborhoods, and choose the tubular neighborhoods and the tame almost complex structure $J$ in the proof of Lemma \ref{lem:Gompf} so that the branch curves, their normal disk fibers, and the curves $E_i$ are all $J$-holomorphic, with $J$ integrable near their intersections. Away from the ramification locus, the lifted almost complex structure preserves $TF_i$. Near $E_i\cap\widehat C_i$ or $E_i\cap\widehat C_i'$, the covering has local coordinates $(u,v)\mapsto (z,w)=(u,v^3)$ with $E_i=\{z=0\}$, and hence $F_i=\{u=0\}$. Thus, $TF_i$ is preserved everywhere by the lifted almost complex structure. Since this structure is tamed by $\omega_Z$, each $F_i$ is $\omega_Z$-symplectic, as desired. 
    
   Since $\pi^*\operatorname{PD}[\widehat C_i] = 3\operatorname{PD}[T_i]$ and $\widehat{C}_i^2=-3$, we have $T_i^2=-1$, and similarly $T_i'^2=-1$. Consider the restriction $\pi:F_i\to E_i$, which is a degree-$3$ branched covering, totally ramified over the two points $E_i\cap \widehat{C}_i$ and $E_i\cap \widehat{C}_i'$. Over the complement of these two points, the restricted covering has monodromy sending positively oriented meridians of the two punctures to $1$ and $2$ in $\mathbb{Z}_3$, respectively. This monodromy is surjective, so $F_i$ is connected. The Riemann--Hurwitz formula implies that $F_i$ is a sphere. Also, since $\pi^*\operatorname{PD}[E_i] = \operatorname{PD}[F_i]$ and $E_i^2=-1$, we see that $F_i^2=-3$. Thus, for $i=1,2,3$, the union $T_i\cup F_i\cup T_i'$ is a configuration of symplectic spheres of type $(-1)$-$(-3)$-$(-1)$.

   Next, let $Y$ be the symplectic manifold obtained from $Z$ by first symplectically blowing down the six spheres $T_i,T_i'$ and then symplectically blowing down the three resulting $(-1)$-spheres obtained from the $F_i$'s. Let $f:Z\to Y$ denote the projection. Note that \[
K_Z= f^*K_Y + \sum_{i=1}^3 \left( \operatorname{PD}[F_i] + 2\operatorname{PD}[T_i] + 2\operatorname{PD}[T_i'] \right).
\]
On the other hand, we have \begin{align*}
    K_Z& =\pi^*K_{\widehat{M}}+2\sum_{i=1}^3(\operatorname{PD}[T_i]+\operatorname{PD}[T_i'])= \pi^*\beta^*K_M+\sum_{i=1}^3(\operatorname{PD}[F_i]+2\operatorname{PD}[T_i]+2\operatorname{PD}[T_i']),
\end{align*}
and hence $f^*K_Y=\pi^*\beta^*K_M$. Since $f^*:H^2(Y;\mathbb{Z})\to H^2(Z;\mathbb{Z})$ is injective and $K_M$ has order two, it follows that $2K_Y=0$, so $Y$ has even intersection form. Also, since \[
e(Y)=e(Z)-9=\left(3e(\widehat{M})-(3-1)\cdot 6\cdot 2 \right)-9=12,
\] 
the formula $0=K_Y^2=2e(Y)+3\sigma(Y)$ implies $\sigma(Y)=-8$. Thus, $Y$ is non-spin by Rohlin's theorem, and hence $K_Y\neq 0$. The adjunction formula shows that $Y$ is minimal. Therefore, Theorem \ref{thm:Li} implies that $Y$ has the integral homology of an Enriques surface. In particular, $Q_Y\cong U\oplus E_8$ and $b_1(Y)=0$. 

Now, consider the fourth $A_2$-configuration $C_4\cup C_4'$ and the $A_1$-sphere $A$. Since their plumbing neighborhoods are simply-connected and disjoint from the branch locus, these configurations lift to three disjoint copies. Thus, there is a $3A_2+3A_1$-configuration in $Z$. Since these spheres are disjoint from the nine blow-down spheres, they descend unchanged to $Y$. It follows that the lattice $3A_2\oplus 3A_1$ embeds in $U\oplus E_8$. We will derive a contradiction by showing that such a lattice embedding cannot exist.

Suppose $S\subset U\oplus E_8$ is a sublattice with $S\cong 3A_2\oplus 3A_1$. The $2$-primary part $(S^*/S)_2$ of its discriminant group is isomorphic to $\mathbb{Z}_2^{\oplus 3}$, and it is easy to see that $q_S(x)$ is $-1/2$, $-1$, or $-3/2 \pmod{2\mathbb{Z}}$ for any nonzero $x\in (S^*/S)_2$. Let $\bar{S}$ denote the primitive closure of $S$ in $U\oplus E_8$. Since the quadratic form $q_S$ vanishes on $\bar{S}/S$,  $(S^*/S)_2 \cap (\bar{S}/S)=0$ and $(\bar{S}/S)$ has an odd order. Furthermore, the bilinear form  $b_S$ on $ (S^*/S)_2 \times \bar{S}/S$ also vanishes, so that $(S^*/S)_2 \subset (\bar{S}/S)^\perp$. Hence, by Lemma \ref{lem:discriminant}(a), the $2$-primary part $(S^*/S)_2$ of the discriminant group is unchanged upon passing to the primitive closure $\bar{S}$ of $S$ in $U\oplus E_8$, i.e., $(\bar{S}^*/\bar{S})_2\cong \mathbb{Z}_2^{\oplus 3}$. 
On the other hand, since $\bar{S}^\perp$ has rank one,  $(\bar{S}^\perp)^*/\bar{S}^\perp$ must be cyclic. But this is a contradiction since Lemma \ref{lem:discriminant}(b) implies that $\bar{S}^*/\bar{S}$ is isomorphic to $(\bar{S}^\perp)^*/\bar{S}^\perp$. 
\end{proof}

\subsection{The non-cyclic case}

In this subsection, we prove the following proposition.

\begin{proposition}\label{prop:non_cyclic_case}
    Let $X$ be a symplectic $\mathbb{Q}$-homology projective plane with only rational double points. If $X$ has at least one non-cyclic rational double point, then $e_{\operatorname{orb}}(X)>0$ and $|\operatorname{Sing}(X)|\leq 5$.
\end{proposition}

This follows from the following two propositions. 

\begin{proposition}\label{prop:D_45A_1}
    There is no symplectic $\mathbb{Q}$-homology projective plane with singularity type $D_4+5A_1$.
\end{proposition}

\begin{proposition}\label{prop:D_4A_23A_1}
    There is no symplectic $\mathbb{Q}$-homology projective plane with singularity type $D_4+A_2+3A_1$.
\end{proposition}

\begin{proof}[Proof of Proposition \ref{prop:non_cyclic_case}]
    Let $n$ denote the number of singular points of $X$, and let $L$ denote the number of the exceptional curves in the minimal symplectic resolution of $X$. Since $X$ contains at least one non-cyclic rational double point, we have $L\geq 4+(n-1)$. Therefore, we have $n\leq 6$ since $L\leq 9$ by Proposition \ref{prop:symp_QHCP2}. 

    If $n\leq 3$, then clearly $e_{\operatorname{orb}}(X)> 0$. If $n\geq 4$, then $X$ must have only one non-cyclic rational double point; otherwise we would have $L\geq 4+4+1+1=10$. Let $n=4$, and suppose that the cyclic rational double points are of types $A_{k_1}$, $A_{k_2}$, $A_{k_3}$. Then $k_1+k_2+k_3\leq 5$, so $\sum_{i=1}^31/(k_i+1)\geq \frac{9}{8}$ by Cauchy--Schwarz. It follows that $e_{\operatorname{orb}}(X)>\frac{1}{8}$.

    If $n=5$ or $n=6$, then $L\leq 9$ implies that the only possible singularity types of $X$ are $D_4+4A_1$, $D_5+4A_1$, $D_4+5A_1$, $D_4+A_2+3A_1$. For the first two cases, we have $e_{\operatorname{orb}}(X)>0$, and the last two cases cannot occur by Propositions \ref{prop:D_45A_1} and \ref{prop:D_4A_23A_1}. In particular, the case $n=6$ does not occur. 
\end{proof}

\begin{proof}[Proof of Proposition \ref{prop:D_45A_1}] 
    The proof is similar to that of Proposition \ref{prop:4A_2A_1}. Suppose $X$ is a symplectic $\mathbb{Q}$-homology projective plane with singularity type $D_4+5A_1$, and let $M\to X$ denote the minimal symplectic resolution. By Proposition \ref{prop:symp_QHCP2}, we have $Q_M\cong U\oplus E_8$, $H_1(M;\mathbb{Z})\cong \mathbb{Z}_2$, and $H_2(M;\mathbb{Z})\cong \mathbb{Z}^{10}\oplus \mathbb{Z}_2$. Thus, $H_2(M;\mathbb{Z}_2)\cong \mathbb{Z}_2^{\oplus 12}$. Since the mod-$2$ intersection pairing is nondegenerate, any totally isotropic subspace $W$ of $H_2(M;\mathbb{Z}_2)$ must satisfy $\dim_{\mathbb{Z}_2} W\leq 6$. 
    
    For the $D_4$-configuration in $M$, let $C$ denote the central sphere and let $L_1$, $L_2$, $L_3$ denote the three outer spheres. Thus, $C^2=L_j^2=-2$, $C\cdot L_j=1$, and the $L_j$'s are pairwise disjoint. Also, let $B_1,\dots,B_5$ denote the (pairwise disjoint) $A_1$-spheres in $M$. We will use the following elementary observations: \begin{enumerate}[label=(\arabic*)]
        \item \label{item_1} For pairwise distinct $S_1,\dots,S_t\in \{L_1,L_2,L_3,B_1,\dots,B_5\}$ with $t\geq 1$, if $\sum_{i=1}^t[S_i]=0\in H_2(M;\mathbb{Z}_2)$, then $4$ divides $t$, and $|\{S_1,\dots,S_t\}\cap \{L_1,L_2,L_3\}|$ is either $0$ or $2$. Hence $t=4$. 
        \item \label{item_2} $\operatorname{dim}_{\mathbb{Z}_2} \operatorname{span}\{[B_1],\dots,[B_5]\}\geq 4$.
    \end{enumerate}
    To prove \ref{item_1}, note that $\sum_{i=1}^t[S_i]=0\in H_2(M;\mathbb{Z}_2)$ implies that $\sum_{i=1}^t [S_i]=2x\in H_2(M;\mathbb{Z})$ for some $x\in H_2(M;\mathbb{Z})$. Since $x^2=Q_M(x,x)$ is even, taking squares gives $4x^2=-2t$. Hence $4$ divides $t$. The second assertion of \ref{item_1} also follows from the observation that \[
    |\{S_1,\dots,S_t\}\cap \{L_1,L_2,L_3\}|=[C]\cdot \sum_{i=1}^t [S_i]=2[C]\cdot x.
    \]
    To prove \ref{item_2}, suppose that the dimension is $\leq 3$. Then we would have two independent relations of the form $[B_{i_1}]+\cdots+[B_{i_4}]=0\in H_2(M;\mathbb{Z}_2)$, and adding the two relations yields a relation of the form $[B_i]+[B_j]=0\in H_2(M;\mathbb{Z}_2)$, which is impossible by \ref{item_1}. 
\smallskip

    Now consider the subspace $W\subset H_2(M;\mathbb{Z}_2)$ generated by the classes of $L_1,L_2,L_3,B_1,\dots,B_5$. Since $W$ is totally isotropic, we have $\dim_{\mathbb{Z}_2}W\leq 6$, so we must have at least two independent linear relations among these classes. By \ref{item_1} and \ref{item_2}, we can find a relation of the form $[L_i]+[L_j]+[B_k]+[B_l]=0$. Therefore, after a suitable relabeling, we may assume that \[
    [L_1]+[L_2]+[B_1]+[B_2]=0\in H_2(M;\mathbb{Z}_2).
    \]
    By Lemma \ref{lem:Baraglia}, there exists a connected double branched cover $\pi:Z\to M$ with branch locus $L_1\cup L_2\cup B_1\cup B_2$. Let $T_1,\dots,T_4$ denote the ramification spheres in $Z$, so that $T_1,T_2,T_3,T_4$ lie over $L_1,L_2,B_1,B_2$, respectively. By Lemma \ref{lem:Gompf}, there exists a symplectic form on $Z$ such that the $T_i$'s are symplectic spheres and \[
    K_Z =\pi^*K_M+\sum_{i=1}^4\operatorname{PD}[T_i]\in H^2(Z;\mathbb{Z}).
    \]
    As in the proof of Proposition \ref{prop:4A_2A_1}, we may assume that $F:=\pi^{-1}(C)\subset Z$ is also symplectic. 

    Note that each $T_i$ is a sphere with $T_i^2=-1$ for $i=1,2,3,4$. Let $f:Z\to Y$ denote the symplectic blow-down of the $T_i$'s, so that \[
    K_Z=f^*K_Y+\sum_{i=1}^4\operatorname{PD}[T_i]\in H^2(Z;\mathbb{Z}).
    \]
    Comparing the two equalities gives $\pi^*K_M=f^*K_Y$. Since $2K_M=0$ (Proposition \ref{prop:symp_QHCP2}) and $f^*$ is injective, we have $2K_Y=0$, so $Y$ has even intersection form. Moreover, $e(Y)=e(Z)-4=(2e(M)-8)-4=12$. Thus, the argument in the proof of Proposition \ref{prop:4A_2A_1} shows that $Y$ has the integral homology of an Enriques surface. In particular, $b_1(Y)=0$. 

    The restriction $\pi|_F:F\to C$ is a connected double branched covering, totally ramified over the two points $C\cap L_1$ and $C\cap L_2$. The Riemann--Hurwitz formula therefore shows that $F$ is a sphere. Since $\pi^*\operatorname{PD}[C]=\operatorname{PD}[F]$, we have $F^2=2C^2=-4$. Moreover, $F$ meets $T_1$ and $T_2$ transversely and positively in one point each and is disjoint from $T_3$ and $T_4$. Consequently, its image $\bar{F}\subset Y$ is a symplectic sphere with $\bar{F}^2=-2$.

    The sphere $L_3$ is disjoint from the branch locus, so it lifts to two disjoint symplectic $(-2)$-spheres in $Z$, and each meets $F$ transversely and positively in one point. Their images, together with $\bar{F}$, form a symplectic $A_3$-configuration in $Y$. Similarly, the remaining spheres $B_3,B_4,B_5$ lift to six disjoint symplectic $(-2)$-spheres in $Z$, disjoint from $F$, the lifts of $L_3$, and the ramification spheres. Their images in $Y$ are therefore disjoint from the preceding $A_3$-configuration. Thus, $Y$ contains a symplectic configuration of type $A_3+6A_1$, whose total rank is $9=b_2^-(Y)$.
    
    Applying Proposition \ref{prop:lattice_constraint} to $Y$ gives six pairwise orthogonal vectors $u_1,\ldots,u_6\in-\mathbb{Z}^8$ satisfying $u_i^2=-2$ and $\mathbf{K}\cdot u_i=2$. As in the proof of Proposition \ref{prop:cyclic_case}, the supports of these vectors are pairwise disjoint and each has cardinality two. Hence there can be at most four such vectors, a contradiction. This proves the proposition.
\end{proof}

\begin{proof}[Proof of Proposition \ref{prop:D_4A_23A_1}] We argue as in the proof of Proposition \ref{prop:D_45A_1}. Suppose $X$ is a symplectic $\mathbb{Q}$-homology projective plane with singularity type $D_4+A_2+3A_1$, and let $M\to X$ be its minimal symplectic resolution. The symplectic manifold $M$ has the integral homology of an Enriques surface, $Q_M\cong U\oplus E_8$, $2K_M=0$, and $H_2(M;\mathbb{Z}_2)\cong \mathbb{Z}_2^{\oplus 12}$.

Let $C$ and $L_1,L_2,L_3$ denote the central and outer spheres of the $D_4$-configuration, respectively. Let $P\cup P'$ denote the $A_2$-configuration, and let $B_1,B_2,B_3$ denote the three $A_1$-spheres. The subspace $V_0\subset H_2(M;\mathbb{Z}_2)$ spanned by the classes of $P$ and $P'$ has intersection matrix \[
\begin{pmatrix}
0&1\\ 1&0
\end{pmatrix},
\]
so it is nondegenerate. Hence $V_0^\perp\subset H_2(M;\mathbb{Z}_2)$ is also nondegenerate and $\dim_{\mathbb{Z}_2}V_0^\perp=10$. The classes of $L_1, L_2, L_3, B_1, B_2, B_3$ span a totally isotropic subspace of $V_0^\perp$, whose dimension is therefore $\leq 5$. Consequently, these six classes satisfy a nonzero relation over $\mathbb{Z}_2$. As in observation \ref{item_1} in the proof of Proposition \ref{prop:D_45A_1}, any such relation involves a number of spheres divisible by $4$ and an even number of outer spheres. Thus it involves exactly $4$ spheres. Since there are only three external $A_1$-spheres, exactly $2$ of the spheres in the relation must be outer spheres. After relabeling, we may therefore assume that $[L_1]+[L_2]+[B_1]+[B_2]=0\in H_2(M;\mathbb{Z}_2)$. Now we perform the same connected double-cover and symplectic blow-down construction as in the proof of Proposition \ref{prop:D_45A_1}. The resulting symplectic manifold $Y$ satisfies $b_1(Y)=0$, $b_2^+(Y)=1$, and $b_2^-(Y)=9$. The $D_4$-configuration gives a symplectic $A_3$-configuration in $Y$. The unused $A_2$-configuration $P\cup P'$ and the remaining $A_1$-sphere $B_3$ lift to two disjoint copies which are unchanged under the blow-downs. Hence $Y$ contains a symplectic configuration of type $A_3+2A_2+2A_1$, whose total rank is $9=b_2^-(Y)$, so that we can apply Proposition \ref{prop:lattice_constraint}. However, we have already shown in the proof of Proposition \ref{prop:A_32A_22A_1} that the corresponding lattice embedding cannot exist. This completes the proof.
\end{proof}

The inequalities in Theorem \ref{thm:main_1} now follow from Proposition \ref{prop:cyclic_case}, Corollary \ref{cor:e_orb_cyclic}, and Proposition \ref{prop:non_cyclic_case}.

We conclude this section with one further exclusion, which will be used in the proof of the classification theorem. The argument is similar to that of Proposition \ref{prop:D_45A_1}. 
\begin{proposition}\label{prop:D_54A_1}
There is no symplectic $\mathbb{Q}$-homology projective plane with singularity type $D_5+4A_1$.
\end{proposition}
\begin{proof}
Suppose that such an $X$ exists. Its minimal symplectic resolution $M$ has the integral homology of an Enriques surface, $Q_M\cong U\oplus E_8$, $2K_M=0$, and $H_2(M;\mathbb{Z}_2)\cong \mathbb{Z}_2^{\oplus 12}$. 

Write the five spheres in the $D_5$-configuration as $C,L_1,L_2,P_1,P_2$, with all self-intersections equal to $-2$ and with the only nonzero intersections between distinct spheres given by \[
C\cdot L_1=C\cdot L_2=C\cdot P_1=P_1\cdot P_2=1.
\]
Let $B_1,\ldots,B_4$ denote the $A_1$-spheres. The subspace $V_0\subset H_2(M;\mathbb{Z}_2)$ spanned by the classes of $P_1$ and $P_2$ is nondegenerate, so $V_0^\perp$ is also nondegenerate and $\dim_{\mathbb{Z}_2}V_0^\perp=10$. The classes of the six pairwise disjoint spheres $L_1,L_2,B_1,B_2,B_3,B_4$ span a totally isotropic subspace of $V_0^\perp$, whose dimension is therefore $\leq 5$. Hence these six classes satisfy a nonzero relation over $\mathbb{Z}_2$. As in observation \ref{item_1} in the proof of Proposition \ref{prop:D_45A_1}, any such relation involves exactly $4$ spheres, and an even number of $L_1, L_2$ occurs in the relation. Therefore, after relabeling the $B_i$'s, the relation is one of \[
[B_1]+[B_2]+[B_3]+[B_4]=0 \quad\text{or}\quad [L_1]+[L_2]+[B_1]+[B_2]=0 \quad\text{in }H_2(M;\mathbb{Z}_2).
\]
In either case, let $B$ be the union of the four spheres in the relation. Take a connected double branched cover $\pi:Z\to M$ with branch locus $B$, and symplectically blow down the four ramification spheres. The resulting symplectic manifold $Y$ satisfies $b_1(Y)=0$, $b_2^+(Y)=1$, and $b_2^-(Y)=9$.  

If $B=B_1\cup\cdots\cup B_4$, then a plumbing neighborhood of the $D_5$-configuration is simply-connected and disjoint from $B$. Its preimage consists of two disjoint copies, which survive the blowdowns unchanged. Hence $Y$ contains a $2D_5$-configuration, forcing $b_2^-(Y)\geq10$, a contradiction.

Suppose instead that $B=L_1\cup L_2\cup B_1\cup B_2$. In this case, we may assume that $F=\pi^{-1}(C)\subset Z$ is symplectic. The image $\bar{F}$ of $F$ in $Y$ is a symplectic $(-2)$-sphere. The $A_2$-chain $P_1\cup P_2$ lifts to two disjoint copies $P_1^+\cup P_2^+$ and $P_1^-\cup P_2^-$. Each $P_1^\pm$ meets $F$ once, whereas each $P_2^\pm$ is disjoint from $F$. After the blow-downs, these curves form the symplectic $A_5$-configuration in $Y$. Considering the unused spheres $B_3,B_4$, we see that $Y$ contains a symplectic configuration of type $A_5+4A_1$.

Proposition \ref{prop:lattice_constraint} now produces four pairwise orthogonal vectors $u_1,\ldots,u_4\in-\mathbb Z^6$ with $u_i^2=-2$ and $\mathbf{K}\cdot u_i=2$. By the support calculation in the proof of Proposition \ref{prop:cyclic_case}, their supports must be pairwise disjoint, each of cardinality two. There can be at most three such vectors in $-\mathbb{Z}^6$, a contradiction. This proves the proposition.
\end{proof}

\section{Proof of Theorem \ref{thm:classification}}\label{sec:classification}

We first reduce the case of five singular points to two possible singularity types.

\begin{proposition}\label{prop:five_singular_points}
Let $X$ be a symplectic $\mathbb{Q}$-homology projective plane with only rational double points. If $|\operatorname{Sing}(X)|=5$, then its singularity type is either $2A_3+3A_1$ or $D_4+4A_1$. 
\end{proposition}
\begin{proof}
Let $L\leq9$ be the number of the exceptional curves in the minimal symplectic resolution.

Suppose first that all singularities are cyclic, and let $k$ be the number of $A_1$-singularities. Proposition \ref{prop:cyclic_case} gives $k \leq3$, while $L\geq k+2(5-k)=10-k$ gives $k\geq1$. If $k=1$, the only possibility is $4A_2+A_1$, which has been excluded in Proposition \ref{prop:4A_2A_1}. If $k=2$, the only possibilities are $3A_2+2A_1$ and $A_3+2A_2+2A_1$. The latter has been excluded in Proposition \ref{prop:A_32A_22A_1}. For the former, $L=8$ and $(9-L)|\det(3A_2\oplus2A_1)|=3^3\cdot2^2$, which is not a square, contrary to Proposition~\ref{prop:determinant_square}. Thus $k=3$.

Apply Proposition~\ref{prop:lattice_constraint} to the resolution. The three $A_1$-singularities give  pairwise orthogonal $u_1,u_2,u_3\in -\mathbb{Z}^6$ with $u_i^2=-2$ and $\mathbf{K}\cdot u_i=2$. As in the proof of Proposition \ref{prop:cyclic_case}, their supports are disjoint and have cardinality two, so we may assume $(u_1,u_2,u_3)=(-e_1+e_2, -e_3-e_4,-e_5-e_6)$. 

Let $A_{k_1}$ and $A_{k_2}$ be the types of the two remaining singularities. Clearly, any vector in $-\mathbb{Z}^6$ orthogonal to $-e_1+e_2, -e_3-e_4,-e_5-e_6$ has even square. It follows that $k_i+1$ is even for $i=1,2$. Since $k_i\geq 2$, we must have $k_i\geq 3$. The bound $L\leq9$ therefore forces $k_1=k_2=3$. Hence the singularity type must be $2A_3+3A_1$.

Suppose now that $X$ has a non-cyclic singularity. There can be only one, by the rank bound. Since it has rank at most five, the only possibilities are $D_4+4A_1$, $D_4+A_2+3A_1$, $D_5+4A_1$. The second and third have been excluded in Propositions \ref{prop:D_4A_23A_1} and \ref{prop:D_54A_1}, respectively. Thus only $D_4+4A_1$ remains. 
\end{proof}

We now consider the case of at most four singular points. Consider the following finite set: \[
\mathcal{T}:=\left\{R_1+\cdots+R_n:1\leq n\leq 4, \, \sum_{i=1}^n \operatorname{rank}R_i\leq 9, \, R_i\in \{A_1,\dots,A_9,D_4,\dots,D_9,E_6,E_7,E_8\}\right\},
\]
where elements that differ only by a permutation of the summands are identified. A direct enumeration gives $|\mathcal T|=127$. This is the same initial list considered in \cite[Section 2.1]{Hwang-Keum-Ohashi-2015}. For completeness, we display the entire list in Table \ref{tab:ade_127}, arranged first by the total rank and then by the number of summands. 

\begingroup
\small
\renewcommand{\arraystretch}{1.4} 
\begin{longtable}{ c | c | >{\raggedright\arraybackslash}p{0.75\textwidth} }
\caption{The 127 types in $\mathcal{T}$ with total rank $L\leq9$ and number of summands $n\leq4$.}
\label{tab:ade_127}\\

\multicolumn{1}{c|}{$L$} & \multicolumn{1}{c|}{$n$} & \multicolumn{1}{c}{Types} \\
\hline \hline
\endfirsthead

\multicolumn{3}{c}{\tablename~\thetable\ (continued)}\\[3pt]
\hline 
\multicolumn{1}{c|}{$L$} & \multicolumn{1}{c|}{$n$} & \multicolumn{1}{c}{Types} \\
\hline \hline
\endhead

\hline
\multicolumn{3}{r}{\footnotesize Continued on the next page}\\
\endfoot

\endlastfoot

1 & 1 & $A_{1}$ \\ \hline
2 & 1 & ${A_{2}}^\dagger$ \\ \hline
2 & 2 & ${2A_{1}}^\dagger$ \\ \hline
3 & 1 & ${A_{3}}^\dagger$ \\ \hline
3 & 2 & $A_{2}+A_{1}$ \\ \hline
3 & 3 & ${3A_{1}}^\dagger$ \\ \hline
4 & 1 & $A_{4}$, ${D_{4}}^\dagger$ \\ \hline
4 & 2 & ${A_{3}+A_{1}}^\dagger$, ${2A_{2}}^\dagger$ \\ \hline
4 & 3 & ${A_{2}+2A_{1}}^\dagger$ \\ \hline
4 & 4 & ${4A_{1}}^\dagger$ \\ \hline
5 & 1 & ${A_{5}}^\dagger$, $D_{5}$ \\ \hline
5 & 2 & ${A_{4}+A_{1}}^\dagger$, ${A_{3}+A_{2}}^\dagger$, ${D_{4}+A_{1}}^\dagger$ \\ \hline
5 & 3 & $A_{3}+2A_{1}$, ${2A_{2}+A_{1}}^\dagger$ \\ \hline
5 & 4 & ${A_{2}+3A_{1}}^\dagger$ \\ \hline
6 & 1 & ${A_{6}}^\dagger$, ${D_{6}}^\dagger$, $E_{6}$ \\ \hline
6 & 2 & $A_{5}+A_{1}$, ${A_{4}+A_{2}}^\dagger$, ${2A_{3}}^\dagger$, ${D_{5}+A_{1}}^\dagger$, ${D_{4}+A_{2}}^\ddagger$ \\ \hline
6 & 3 & ${A_{4}+2A_{1}}^\dagger$, ${A_{3}+A_{2}+A_{1}}^\dagger$, $3A_{2}$, ${D_{4}+2A_{1}}^\dagger$ \\ \hline
6 & 4 & ${A_{3}+3A_{1}}^\dagger$, ${2A_{2}+2A_{1}}^\dagger$ \\ \hline
7 & 1 & $A_{7}$, ${D_{7}}^\dagger$, $E_{7}$ \\ \hline
7 & 2 & ${A_{6}+A_{1}}^\dagger$, $A_{5}+A_{2}$, ${A_{4}+A_{3}}^\dagger$, $D_{6}+A_{1}$, ${D_{5}+A_{2}}^\dagger$, ${D_{4}+A_{3}}^\dagger$, ${E_{6}+A_{1}}^\dagger$ \\ \hline
7 & 3 & ${A_{5}+2A_{1}}^\dagger$, ${A_{4}+A_{2}+A_{1}}^\dagger$, $2A_{3}+A_{1}$, ${A_{3}+2A_{2}}^\dagger$, ${D_{5}+2A_{1}}^\dagger$, ${D_{4}+A_{2}+A_{1}}^\dagger$ \\ \hline
7 & 4 & ${A_{4}+3A_{1}}^\dagger$, ${A_{3}+A_{2}+2A_{1}}^\dagger$, ${3A_{2}+A_{1}}^\dagger$, $D_{4}+3A_{1}$ \\ \hline
8 & 1 & $A_{8}$, $D_{8}$, $E_{8}$ \\ \hline
8 & 2 & $A_{7}+A_{1}$, ${A_{6}+A_{2}}^\dagger$, ${A_{5}+A_{3}}^\dagger$, $2A_{4}$, ${D_{7}+A_{1}}^\dagger$, ${D_{6}+A_{2}}^\dagger$, $D_{5}+A_{3}$, ${D_{4}+A_{4}}^\dagger$, $2D_{4}$, $E_{7}+A_{1}$, $E_{6}+A_{2}$ \\ \hline
8 & 3 & ${A_{6}+2A_{1}}^\dagger$, $A_{5}+A_{2}+A_{1}$, ${A_{4}+A_{3}+A_{1}}^\dagger$, ${A_{4}+2A_{2}}^\dagger$, ${2A_{3}+A_{2}}^\dagger$, $D_{6}+2A_{1}$, ${D_{5}+A_{2}+A_{1}}^\dagger$, ${D_{4}+A_{3}+A_{1}}^\dagger$, ${D_{4}+2A_{2}}^\ddagger$, ${E_{6}+2A_{1}}^\dagger$ \\ \hline
8 & 4 & ${A_{5}+3A_{1}}^\dagger$, ${A_{4}+A_{2}+2A_{1}}^\dagger$, $2A_{3}+2A_{1}$, ${A_{3}+2A_{2}+A_{1}}^\dagger$, $4A_{2}$, ${D_{5}+3A_{1}}^\dagger$, ${D_{4}+A_{2}+2A_{1}}^\dagger$ \\ \hline
9 & 1 & $A_{9}$, $D_{9}$ \\ \hline
9 & 2 & $A_{8}+A_{1}$, $A_{7}+A_{2}$, ${A_{6}+A_{3}}^\ddagger$, $A_{5}+A_{4}$, $D_{8}+A_{1}$, ${D_{7}+A_{2}}^\ddagger$, $D_{6}+A_{3}$, $D_{5}+A_{4}$, $D_{5}+D_{4}$, ${D_{4}+A_{5}}^\ddagger$, $E_{8}+A_{1}$, $E_{7}+A_{2}$, $E_{6}+A_{3}$ \\ \hline
9 & 3 & $A_{7}+2A_{1}$, $A_{6}+A_{2}+A_{1}$, $A_{5}+A_{3}+A_{1}$, $A_{5}+2A_{2}$, $2A_{4}+A_{1}$, ${A_{4}+A_{3}+A_{2}}^\ddagger$, $3A_{3}$, $D_{7}+2A_{1}$, $D_{6}+A_{2}+A_{1}$, $D_{5}+A_{3}+A_{1}$, ${D_{5}+2A_{2}}^\ddagger$, ${D_{4}+A_{4}+A_{1}}^\ddagger$, $2D_{4}+A_{1}$, ${D_{4}+A_{3}+A_{2}}^\ddagger$, $E_{7}+2A_{1}$, $E_{6}+A_{2}+A_{1}$ \\ \hline
9 & 4 & ${A_{6}+3A_{1}}^\ddagger$, $A_{5}+A_{2}+2A_{1}$, $A_{4}+A_{3}+2A_{1}$, ${A_{4}+2A_{2}+A_{1}}^\ddagger$, $2A_{3}+A_{2}+A_{1}$, $A_{3}+3A_{2}$, $D_{6}+3A_{1}$, ${D_{5}+A_{2}+2A_{1}}^\ddagger$, $D_{4}+A_{3}+2A_{1}$, ${D_{4}+2A_{2}+A_{1}}^\ddagger$, ${E_{6}+3A_{1}}^\ddagger$ 
\end{longtable}
\endgroup 

\begin{proposition}\label{prop:at_most_four_singular_points}
Let $X$ be a singular symplectic $\mathbb{Q}$-homology projective plane with only rational double points. If $|\operatorname{Sing}(X)|\leq 4$, then its singularity type is one of the $57$ unmarked types in \mbox{Table~\textup{\ref{tab:ade_127}}}.
\end{proposition} 
\begin{proof}
Let $M\to X$ be the minimal symplectic resolution, and let $R\subset Q_M$ be the sublattice generated by the classes of the exceptional curves. Put $L:=\operatorname{rank}R$. Proposition \ref{prop:symp_QHCP2} gives $L\leq9$, so the singularity type belongs to $\mathcal{T}$.

If $L<9$, Proposition \ref{prop:determinant_square} requires $(9-L)|\det R|$ to be a square. Using \[
|\det A_k|=k+1,\quad |\det D_k|=4,\quad |\det E_6|=3,\quad |\det E_7|=2,\quad |\det E_8|=1,
\]
we find that precisely the $56$ types marked with a dagger fail this condition; they agree with the types in \cite[Table 2]{Hwang-Keum-Ohashi-2015}. These types are therefore excluded. 

For $L=9$, Proposition \ref{prop:symp_QHCP2} gives $Q_M\cong U\oplus E_8$. The twelve rank-$9$ types marked with a double dagger do not embed in this lattice, by \cite[Lemma 2.4(1)]{Hwang-Keum-Ohashi-2015}. For $2\leq L<9$, Proposition \ref{prop:determinant_square} gives $Q_M\cong \langle+1\rangle\oplus L\langle-1\rangle$. The two remaining double-dagger types, $D_4+A_2$ and $D_4+2A_2$, are excluded by \cite[Lemma 2.4(2)]{Hwang-Keum-Ohashi-2015}. Thus, a total of $70=56+14$ types are excluded from the initial list of $127$ types, leaving precisely the $57$ unmarked types in the table.
\end{proof}

Propositions \ref{prop:five_singular_points} and \ref{prop:at_most_four_singular_points} complete the restriction on singularity types.  

\begin{proof}[Proof of Theorem \ref{thm:classification}]
We have $|\operatorname{Sing}(X)|\leq5$ by Theorem \ref{thm:main_1}. If $X$ has at most four singular points, Proposition \ref{prop:at_most_four_singular_points} shows that its type is one of the $57$ unmarked types in Table \ref{tab:ade_127}. If $X$ has five singular points, Proposition \ref{prop:five_singular_points} leaves only $2A_3+3A_1$ and $D_4+4A_1$. These give $57+2=59$ distinct candidate types, precisely those in Table \ref{tab:classification}.

Among the $57$ types with at most four singular points, $27$ have rank $<9$ and $30$ have rank $9$. Together with the type $2A_3+3A_1$, these give $58$ algebraically realizable types. The $27$ types of rank $<9$ are realized by Gorenstein log del Pezzo surfaces of Picard number one; see \cite[Section 1]{Hwang-Keum-Ohashi-2015} and the references therein. The $30$ rank-$9$ types, together with $2A_3+3A_1$, are realized by contracting configurations of rational curves on Enriques surfaces; see \cite[Proposition 8.2]{Schutt-2023}. Equipped with K\"ahler orbifold forms, these projective surfaces give symplectic $\mathbb{Q}$-homology projective planes with the prescribed singularity types. Thus all $58$ types other than $D_4+4A_1$ are realized, proving the assertion.
\end{proof}

Now we prove Corollaries \ref{cor:odd_homology_classification} and \ref{cor:smooth_locus_homology}.

\begin{proof}[Proof of Corollary \ref{cor:odd_homology_classification}]
Let $\rho:M\to X$ be the minimal symplectic resolution, and let $L$ be the total rank of its exceptional configuration. Replacing a simply-connected plumbing by the cone on its boundary does not change the fundamental group. Hence $\pi_1(M)\cong \pi_1(X)$. 

If $L=9$, Proposition \ref{prop:symp_QHCP2} implies that $M$ has the integral homology of an Enriques surface. In particular, $H_1(X;\mathbb{Z}_2)\cong\mathbb {Z}_2$, contrary to the hypothesis. The type $D_4+4A_1$ is excluded by Proposition \ref{prop:exceptional_intro}. Theorem \ref{thm:classification} therefore leaves precisely the stated $27$ types.

Each of these types is realized by a Gorenstein log del Pezzo surface of Picard number one; see \cite[Section 1]{Hwang-Keum-Ohashi-2015} and the references therein. Its minimal resolution is a rational surface and hence simply-connected. Contracting the exceptional plumbings preserves the fundamental group, so the underlying space of the resulting projective surface is simply-connected. Equipped with an orbifold K\"ahler form, it gives the required realization.
\end{proof}

\begin{proof}[Proof of Corollary \ref{cor:smooth_locus_homology}]
Let $M\to X$ be the minimal symplectic resolution. Let $P$ be a disjoint union of closed regular neighborhoods of the exceptional configurations, and put $M_0=M-\operatorname{int}P$. Let $R\subset Q_M$ denote the exceptional lattice and $L:=\operatorname{rank}R$. Since $M_0$ is homotopy equivalent to $X^\circ$, the exact sequence of the pair $(M,M_0)$ can be written as \[
H_2(M;\mathbb Z)\longrightarrow R^*
\longrightarrow H_1(X^\circ;\mathbb Z)
\longrightarrow H_1(M;\mathbb Z)\longrightarrow0.
\]
The hypothesis implies $H_1(M;\mathbb Z)=0$ and that the intersection restriction map $H_2(M;\mathbb Z)\to R^*$ is surjective. Since $Q_M$ is unimodular, this surjectivity implies that $R$ is primitive in $Q_M$. The case $L=9$ is impossible by the homology conclusion of Proposition \ref{prop:symp_QHCP2}.

If $L=1$, the only type is $A_1$. Suppose $2\leq L\leq8$, and let $k\in H_2(M;\mathbb{Z})/\text{tor}$ be the Poincar\'e dual of the canonical class modulo torsion. Then $k^2=9-L\in\{1,\dots,7\}$ and $k\perp R$. The unimodular lattice $Q_M$ is odd, since its signature $1-L$ is not divisible by $8$. As $k$ is characteristic, it cannot be an even multiple of an integral vector. Nor can $k=mv$ for an odd integer $m\geq3$, since $v^2$ would be a positive integer and hence $k^2=m^2v^2\geq9$. Hence $k$ is primitive.

Both $R$ and $k^\perp$ have rank $L$, and $R$ is primitive, so $R=k^\perp$. Unimodularity of $Q_M$ now gives $|\det R|=k^2=9-L$. Among the rank-$L$ entries of Table \ref{tab:classification}, with $2\leq L\leq8$, this equality holds precisely for
\[
A_2+A_1,\quad A_4,\quad D_5,\quad E_6,\quad E_7,\quad E_8.
\]
Together with $A_1$, these are the seven claimed types. All seven types have projective realizations with a simply-connected smooth locus, by the classification of Gorenstein log del Pezzo surfaces with simply-connected smooth locus; see \cite[Theorem 1.6]{Jo-Park-Park-2024-02} and the references therein. 
\end{proof}

\section{The exceptional type \texorpdfstring{$D_4+4A_1$}{D4+4A1}}\label{sec:exceptional_type}

In this section, we study the remaining type $D_4+4A_1$ in Theorem \ref{thm:classification}. We establish the following restrictions on a hypothetical symplectic realization.

\begin{proposition}\label{prop:exceptional_intro}
Let $(X,\omega)$ be a symplectic $\mathbb{Q}$-homology projective plane with singularity type $D_4+4A_1$, and let $\rho:(M,\omega_M)\to(X,\omega)$ be its minimal symplectic resolution. Then $H_1(X;\mathbb{Z}_2)\neq0$ and
\[
|H_1(X^\circ;\mathbb{Z})|
=8\cdot|H_1(X;\mathbb{Z})|\geq16.
\]
Moreover, $\pi_1(X^\circ)$ is infinite, $(M,\omega_M)$ is symplectically minimal, $K_M^2=1$, and $\kappa(M,\omega_M)=2$.
\end{proposition}

We prove Proposition \ref{prop:exceptional_intro} in the following subsections. Fix $X$ and $\rho:(M,\omega_M)\to(X,\omega)$ as in the statement. Let $P$ be the disjoint union of plumbing neighborhoods of the exceptional configurations, and put $M_0:=M-\operatorname{int}P$. Note that $M_0$ is homotopy equivalent to the smooth locus $X^\circ$. Replacing each simply-connected plumbing by the cone on its boundary does not change the fundamental group, so $\pi_1(M)\cong\pi_1(X)$. Proposition \ref{prop:symp_QHCP2} gives
\[
b_1(M)=0,\quad b_2^+(M)=1,\quad b_2^-(M)=8,\quad
e(M)=11,\quad \sigma(M)=-7,\quad K_M^2=1.
\]

\subsection{The exceptional lattice} Write $\Lambda=H_2(M;\mathbb Z)/\text{tor}$ and let
$k\in\Lambda$ denote the Poincar\'e dual of $K_M$ modulo torsion. The vector $k$ is characteristic, has square one, and is orthogonal to the exceptional lattice $R\cong D_4\oplus 4A_1$. Since $k^2=1$, there is an orthogonal decomposition $\Lambda=\mathbb Z\cdot k\oplus E$ with $E:=k^\perp \cong E_8$. Indeed, $E$ is negative definite, even, unimodular, and of rank $8$. Moreover, $|\det R|=64$, $[E:R]=8$, and $R^*/R\cong \mathbb{Z}_2^{\oplus 6}$. The subgroup $H:=E/R\subset R^*/R$ is isotropic with respect to the discriminant quadratic form $q_R$ and has dimension $3$ over $\mathbb{Z}_2$.

Let $B_1,B_2,B_3,B_4\subset M$ denote the four $A_1$-spheres. The composition \[
H=E/R \lhook\joinrel\longrightarrow R^*/R \xrightarrow{\operatorname{projection}} D_4^*/D_4\cong \mathbb{Z}_2^{\oplus 2}
\] 
has nonzero kernel. A class supported on $w$ of the four $A_1$ coordinates has discriminant square $-w/2$ modulo $2\mathbb{Z}$. The only nonzero isotropic class supported on these coordinates is therefore the class supported on all four. Consequently, \begin{equation}\label{eq:exceptional_half_sum}
\frac{[B_1]+[B_2]+[B_3]+[B_4]}2\in E\subset\Lambda.
\end{equation}

\subsection{Nontrivial mod-2 homology} Suppose, for a contradiction, that $H_1(M;\mathbb{Z}_2)=0$. Then $H_1(M;\mathbb{Z})$ is finite of odd order, and therefore $\operatorname{tor}H_2(M;\mathbb{Z})$ also has odd order. Thus, \eqref{eq:exceptional_half_sum} lifts to a relation \[
[B_1]+[B_2]+[B_3]+[B_4]=2d_0\in H_2(M;\mathbb{Z}).
\]
Any three of the four classes are independent in $H_2(M;\mathbb{Z}_2)$. Indeed, a nonzero relation supported on a proper subset $I\subset\{1,2,3,4\}$ would give \[
\sum_{i\in I}\,[B_i]=2x\in H_2(M;\mathbb{Z}).
\]
Pairing with $K_M$ gives $\langle K_M, x\rangle=0$, so the image of $x$ in $\Lambda$ lies in the even lattice $E=k^\perp$. Taking squares yields $4x^2=-2|I|$, which implies that $4$ divides $|I|$, a contradiction.

By Lemma \ref{lem:Baraglia}, there is a connected double cover $\pi:Z\to M$ branched over the four spheres, with $b_1(Z)=0$. Lemma \ref{lem:Gompf} equips $Z$ with a symplectic form for which the ramification spheres $T_i=\pi^{-1}(B_i)$ are symplectic and satisfy $T_i^2=-1$. Let $f:Z\to Y$ denote their symplectic blow-down. The canonical classes satisfy \[
K_Z=\pi^*K_M+\sum_{i=1}^4\operatorname{PD}[T_i]
=f^*K_Y+\sum_{i=1}^4\operatorname{PD}[T_i],
\]
so $f^*K_Y=\pi^*K_M$ and $K_Y^2=2$. Since $e(Y)=10$ and $b_1(Y)=0$, it follows from $K_Y^2=2e(Y)+3\sigma(Y)$ that $\sigma(Y)=-6$ and $b_2^-(Y)=7$. The $D_4$-configuration in $M$ is disjoint from the branch locus, so it lifts to two disjoint copies in $Z$. Both copies survive the blow-downs and give a rank-$8$ negative definite sublattice of $Q_Y$, a contradiction. We conclude that $H_1(X;\mathbb{Z}_2)\cong H_1(M;\mathbb{Z}_2)\neq 0$.

\subsection{Homology of the smooth locus}
The homology exact sequence of the pair $(M,M_0)$ can be written as \[
H_2(M;\mathbb{Z})\longrightarrow R^* \longrightarrow H_1(X^\circ;\mathbb{Z}) \longrightarrow H_1(M;\mathbb{Z})\longrightarrow0.
\]
The image of the first map is precisely $E\subset R^*$. Since $R^*/E\cong \mathbb{Z}_2^{\oplus 3}$, we obtain the short exact sequence \begin{equation}\label{eq:exceptional_homology_sequence}
0\longrightarrow \mathbb{Z}_2^{\oplus 3} \longrightarrow H_1(X^\circ;\mathbb{Z}) \longrightarrow H_1(X;\mathbb {Z})\longrightarrow0.
\end{equation}
All groups here are finite, and the last has even order. This proves \[
|H_1(X^\circ;\mathbb {Z})|=8\cdot |H_1(X;\mathbb{Z})|\geq 16.
\]

\subsection{Minimality and Kodaira dimension}
Since $b_1(M)=0$ and $H_1(M;\mathbb{Z}_2)\neq0$, $M$ is neither rational nor ruled. A symplectic minimal model obtained by $r$ blow-downs has canonical square $1+r>0$, and therefore has Kodaira dimension two. Hence $\kappa(M,\omega_M)=2$. Note that $\kappa$ depends only on the oriented diffeomorphism type of $M$. For more on symplectic Kodaira dimension, see \cite{Li-2019}. 

We next prove global minimality. The construction in the proof of Theorem \ref{thm:symplectic_resolution} gives a smooth family $\omega_\epsilon$ on $M$ as in Remark \ref{rmk:resolution_family}, satisfying \eqref{eq:resolution_family}. In particular, \[
[\omega_\epsilon]=\alpha+\epsilon\eta\longrightarrow\alpha, \quad \alpha:=\rho^*[\omega].
\]
Note that $\alpha^2>0$ and $\alpha|_R=0$. The canonical class is constant along the family $\{\omega_\epsilon\}$. Write $k_0\in\Lambda$ for the Poincar\'e dual of its common canonical class modulo torsion, and put $E_0=k_0^\perp$. As above, $k_0^2=1$, $E_0\cong E_8$, and $R\otimes\mathbb R=E_0\otimes\mathbb R$. Each $(M,\omega_\epsilon)$ has Kodaira dimension two, so $\langle [\omega_\epsilon],k_0\rangle = K_{\omega_\epsilon}\cdot [\omega_\epsilon]>0$ (this positivity does not require minimality of $(M,\omega_\epsilon)$), which implies $\langle \alpha,k_0\rangle\geq 0$. We cannot have $\langle \alpha,k_0\rangle=0$ since $\alpha^2>0$ and $\alpha|_{E_0}=0$. Hence $\langle \alpha,k_0\rangle>0$.

If $e\in H_2(M;\mathbb{Z})/\text{tor}$ is represented by an $\omega_\epsilon$-symplectic $(-1)$-sphere, the adjunction formula gives $k_0\cdot e=-1$. It follows that $e=-k_0+v$ for some $v\in E_0$ with $v^2=-2$. There are only finitely many such $v$, and for every one of them, we have $\langle [\omega_\epsilon],-k_0+v\rangle\to \langle\alpha,-k_0+v\rangle=-\langle \alpha,k_0\rangle<0$. Consequently, for all sufficiently small $\epsilon$, no symplectic $(-1)$-sphere exists. Recall that symplectic and smooth minimality are equivalent for closed symplectic $4$-manifolds \cite[Proposition 2.1]{Li-2006}. It follows that $M$ is smoothly minimal. In particular, the original resolution $(M,\omega_M)$ is symplectically minimal as well.

\subsection{The fundamental group of the smooth locus}
Put $G=H_1(X^\circ;\mathbb{Z})$ and $d=|G|$. Take the connected cover of $X^\circ$ defined by the abelianization homomorphism $\pi_1(X^\circ)\to G$. The kernel of the homomorphism \[
R^*/R\cong H_1(\partial P;\mathbb{Z})\longrightarrow H_1(X^\circ;\mathbb{Z})=G
\]
is $H=E/R$. Every nonzero class in $D_4^*/D_4$ has discriminant quadratic value $-1\pmod{2\mathbb Z}$, while the unique nonzero class in each $A_1^*/A_1$ has value $-1/2\pmod{2\mathbb Z}$. Since $H$ is isotropic with respect to $q_R$, its intersection with each summand of $R^*/R\cong (D_4^*/D_4)\oplus 4(A_1^*/A_1)$ is trivial. Thus the restriction of this homomorphism to each summand is injective. Identifying each summand with the abelianization of the corresponding isotropy group $\Gamma_p$, we conclude that the composition $\Gamma_p\to\Gamma_p^{\mathrm{ab}}\to G$ has kernel $[\Gamma_p,\Gamma_p]$.

Complete the cover using the corresponding quotient-ball charts, equipped with the lifted symplectic forms. The lifts of each $A_1$-singularity are smooth, and there are $d/2$ points above it. The isotropy group at the $D_4$-singularity is the quaternion group $Q_8$, with $[Q_8,Q_8]=\{\pm1\}$ and $Q_8^{\mathrm{ab}}\cong\mathbb{Z}_2^{\oplus2}$. Thus there are $d/4$ points above this singularity, all of type $A_1$. Consequently, the completed covering $\widehat{X}\to X$ is a connected symplectic orbifold covering of degree $d$ with exactly $d/4$ singular points, all of type $A_1$.

Let $\widetilde{X}\to\widehat{X}$ be the minimal symplectic resolution. Since $e_{\mathrm{orb}}(X)=1/8$, multiplicativity of the orbifold Euler characteristic gives $e_{\mathrm{orb}}(\widehat{X})=d/8$. Hence $e(\widehat{X})=d/4$ and $e(\widetilde{X})=d/2$.

The orbifold canonical class pulls back under an orbifold covering, and resolving rational double points introduces no exceptional correction terms in the canonical class. Therefore, \[
K_{\widetilde{X}}^2=K_{\widehat{X}}^2=dK_X^2=dK_M^2=d,
\]
and this implies $\sigma(\widetilde{X})=0$. Writing $t=b_1(\widetilde{X})$, we obtain $b_2^-(\widetilde{X})=b_2^+(\widetilde{X})=d/4+t-1$. The $d/4$ disjoint exceptional $(-2)$-spheres in $\widetilde{X}$ force $b_2^-(\widetilde{X})\geq d/4$, hence $t\geq1$. If $\pi_1(X^\circ)$ were finite, then $\pi_1(\widetilde{X})$, being a quotient of a subgroup of $\pi_1(X^\circ)$, would also be finite, contradicting $b_1(\widetilde{X})\geq1$. We conclude that $\pi_1(X^\circ)$ is infinite. This completes the proof of Proposition \ref{prop:exceptional_intro}.

\begin{remark}\label{rmk:smooth_exceptional_type} There exists an oriented smooth orbifold with singularity type $D_4+4A_1$ whose underlying space is simply-connected and has the rational homology of $\mathbb{CP}^2$. By Proposition \ref{prop:exceptional_intro}, this orbifold admits no symplectic structure compatible with its orientation and orbifold structure.

Start with a Gorenstein log del Pezzo surface $S$ of Picard number one and singularity type $2D_4$. Let $P(D_4)$ denote the negative definite $D_4$-plumbing. In the notation of \cite[Sections 2.3--2.4]{Baraglia-2026}, the disk bundle $R(2)$ over $\mathbb{RP}^2$ has oriented total space, and $\partial R(2)\cong \partial P(D_4)$ as oriented $3$-manifolds. Since $R(2)$ deformation retracts onto $\mathbb{RP}^2$, it is a rational homology $4$-ball. Replace a cone neighborhood of one singular point of $S$ by $R(2)$. The resulting orbifold $S'$ has the rational homology of $\mathbb{CP}^2$ and one singularity of type $D_4$.

We claim that $S'$ is simply-connected. Let $\tilde{S}\to S$ be the minimal resolution, and let $P_1,P_2\subset \tilde{S}$ be disjoint, closed, regular neighborhoods of the two $D_4$-configurations. Put $M_0=\tilde{S}-\operatorname{int}(P_1\cup P_2)$ and $C=\tilde{S}-\operatorname{int}P_1$. The surface $\tilde{S}$ is rational and hence simply-connected. Note that the lattice $D:=H_2(P_1;\mathbb{Z})\cong D_4$ is primitive in $H_2(\tilde{S};\mathbb{Z})$. Indeed, its primitive closure $\overline{D}$ is even, since $D$ is orthogonal to the canonical class by the adjunction formula. Since the discriminant quadratic form of $D_4$ has no nonzero isotropic element, we have $\overline{D}=D$ by Lemma~\ref{lem:discriminant}(a). The exact sequence of the pair $(\tilde{S},C)$ therefore gives $H_1(C;\mathbb{Z})=0$.

By \cite[Lemma 6]{Miyanishi-Zhang-1988}, we have $\pi_1(M_0)\cong\pi_1(S^\circ)\cong \mathbb{Z}_2^{\oplus 2}$. Since $C$ is obtained from $M_0$ by attaching the simply-connected plumbing $P_2$, its fundamental group is a quotient of $\pi_1(M_0)$. Thus, $\pi_1(C)$ is abelian, and $H_1(C;\mathbb Z)=0$ implies $\pi_1(C)=1$. Contracting the $D_4$-configuration in $P_2\subset C$ does not change the fundamental group. Finally, the homomorphism $\pi_1(\partial R(2))\to\pi_1(R(2))$ is surjective since the circle bundle $\partial R(2)\to\mathbb{RP}^2$ induces a surjection on $\pi_1$. Van Kampen's theorem now shows that $\pi_1(S')=1$.

Recall the oriented smooth orbifold $Q=S^4/\langle\tau\rangle$ from the introduction, which is simply-connected and has the rational homology of $S^4$. Its two singularities are of type $A_1$. Now the space $X:=S'\#2Q$, where the connected sums are taken in the smooth loci, is an oriented, smooth, $\mathbb{Q}$-homology projective plane with singularity type $D_4+4A_1$ and $\pi_1(X)=1$.
\end{remark}

\begin{remark}
    The restrictions in Proposition \ref{prop:exceptional_intro} provide further motivation for expecting that $D_4+4A_1$ is not symplectically realizable. A related recent result of Chen--Shin \cite[Theorem 1.1]{Chen-Shin-2026} states that a numerical Godeaux surface contains at most six disjoint $(-2)$-curves. For a symplectic realization of $D_4+4A_1$, the three outer spheres of the $D_4$-configuration, together with the four $A_1$-spheres, would form seven pairwise disjoint symplectic $(-2)$-spheres in the minimal symplectic resolution, which has $b_1=0$, $b_2^+=1$, $K^2=1$, and symplectic Kodaira dimension two. Thus, a symplectic analog of their bound would exclude the type $D_4+4A_1$.
\end{remark}

\appendix
\section{Proof of Proposition \ref{prop:lattice_constraint}}\label{sec:appendix}

In this appendix, we prove Proposition \ref{prop:lattice_constraint} by adapting the $\operatorname{spin}^c$-gluing argument of \cite{Jo-Park-Park-2026}. We first fix our conventions for lens spaces and plumbings and recall the boundary restriction formulas needed in the proof, following \cite[Section 3]{Jo-Park-Park-2026}. For background on $\operatorname{spin}^c$ structures, see \cite[Section 2.1]{Jo-Park-Park-2026}.

Throughout the preliminary discussion, let $p>q>0$ be coprime integers.
\begin{itemize}
    \item The lens space $L(p,q)$ is defined to be the quotient \[
L(p,q):=\frac{S^3}{(z,w)\sim (\zeta z,\zeta^qw)}  \quad (\zeta=e^{2\pi i /p}),
\]
where we view $S^3$ as the unit sphere in $\mathbb{C}^2$. We orient $L(p,q)$ by the quotient orientation induced from the standard orientation of $S^3$. 
\item Suppose that the Hirzebruch--Jung continued fraction of $p/q$ is given by \[ \frac{p}{q} = [a_1, \dots, a_s] := a_1 - \frac{1}{a_2 - \dfrac{1}{\cdots - \dfrac{1}{a_s}}} \quad (a_i \geq 2). \]
Then $X(p,q)$ denotes the negative definite plumbed $4$-manifold with $\partial X(p,q)=L(p,q)$, constructed according to the linear graph shown in the following figure. 
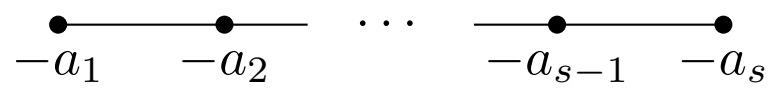
\begin{figure}[h]
\centering
\begin{tikzpicture}[scale=1.1]
\draw (-2,0) node[circle, fill, inner sep=1.2pt, black]{};
\draw (-1,0) node[circle, fill, inner sep=1.2pt, black]{};
\draw (1,0) node[circle, fill, inner sep=1.2pt, black]{};
\draw (2,0) node[circle, fill, inner sep=1.2pt, black]{};
\draw (-2,0) node[below]{$-a_1$};
\draw (-1,0) node[below]{$-a_2$};
\draw (1,0) node[below]{$-a_{s-1}$};
\draw (2,0) node[below]{$-a_s$};
\draw (0,0) node{$\cdots$};
\draw (-2,0)--(-1,0) (-1,0)--(-0.5,0) (0.5,0)--(1,0)  (1,0)--(2,0) ;
\end{tikzpicture}
\end{figure}

\noindent The corresponding homology classes $v_1,\dots,v_s$ form a basis for $H_2(X(p,q);\mathbb{Z})$, called the \emph{standard basis}.
\item There is an identification of the set $\operatorname{Spin}^c(L(p,q))$ of $\operatorname{spin}^c$ structures on $L(p,q)$ with $\mathbb{Z}_p$, provided by Ozsv\'ath and Szab\'o \cite[Section 4.1]{Ozsvath-Szabo-2003}. We denote by $\mathfrak{s}_{p,q}(i) \in \operatorname{Spin}^c(L(p,q))$ the $\operatorname{spin}^c$ structure corresponding to $i\in \mathbb{Z}_p$ under this identification. 
\item For $\epsilon\in \{\pm 1\}$, the Fukumoto-Furuta-Ue $\sigma$-function \cite{Fukumoto-Furuta-Ue-2001} is defined by \[
  \sigma(q,p,\epsilon)=\frac{1}{p}\sum_{k=1}^{p-1} \left(\cot \frac{\pi k}{p}\cot \frac{\pi k q}{p} +2\epsilon^k \csc \frac{\pi k }{p}\csc \frac{\pi k q}{p} \right).
\]
For even $p$, we have $\sigma (q,p,+1)=\sigma (q,p,-1)$ precisely when the two spin structures of the lens space $L(p,q)$ have the same Heegaard Floer $d$-invariant. It is easy to see directly that $\sigma (q,p,+1)\neq \sigma (q,p,-1)$ if $p$ is even and $q\in \{1,p-1\}$. 
\end{itemize}

From now on, we assume that either $p$ is odd or $\sigma(q,p,+1)\neq \sigma(q,p,-1)$.  
\begin{itemize}
    \item Let $f_{p,q}\colon L(p,q)\to L(p,p-q)$ be the orientation-reversing diffeomorphism defined by $f_{p,q}[z_1,z_2]=[\bar{z}_1,z_2]$, where $\bar{z}_1$ denotes the complex conjugate of $z_1$. Then \begin{equation}\label{eq:pullback}
        f_{p,q}^*\mathfrak{s}_{p,p-q}(i)=\mathfrak{s}_{p,q}(p-1-i)\in \operatorname{Spin}^c(L(p,q))
    \end{equation}
    for all $i\in \mathbb{Z}_p$ \cite[Proposition 3.8]{Jo-Park-Park-2026}. 
    \item For the standard basis $\{v_1,\dots,v_s\}$ for $H_2(X(p,q);\mathbb{Z})$, let $\{v_1^*,\dots,v_s^*\}$ denote the dual basis for $H^2(X(p,q);\mathbb{Z})$ so that $\langle v_i^*,v_j\rangle=\delta_{ij}$. Consider the cohomology classes \[
    K^0_{p,q}:=\sum_{i=1}^s (a_i-2)v_i^*, \quad K_{p,q}:=2v_s^*+K^0_{p,q}=2v_s^*+\sum_{i=1}^s (a_i-2)v_i^* \in H^2(X(p,q);\mathbb{Z}). 
    \]
    Then there exist unique $\operatorname{spin}^c$ structures $\mathfrak{s}_{K_{p,q}^0}$ and $\mathfrak{s}_{K_{p,q}}$ on $X(p,q)$ such that $c_1(\mathfrak{s}_{K_{p,q}^0})=K^0_{p,q}$ and $c_1(\mathfrak{s}_{K_{p,q}})=K_{p,q}$. Under the restriction map $\operatorname{Spin}^c(X(p,q))\to \operatorname{Spin}^c(L(p,q))$, we have \begin{equation}\label{eq:restriction}
         \mathfrak{s}_{K_{p,q}^0}|_{L(p,q)}=\mathfrak{s}_{p,q}(p-1), \quad \mathfrak{s}_{K_{p,q}}|_{L(p,q)}=\mathfrak{s}_{p,q}(0)\in \operatorname{Spin}^c(L(p,q));
    \end{equation}
    see the proof of \cite[Proposition 3.5]{Jo-Park-Park-2026}.
\end{itemize}

\begin{proof}[Proof of Proposition~\ref{prop:lattice_constraint}]
Put $p_i:=k_i+1$ and $N:=\sum_{i=1}^nk_i=b_2^-(M)$. For each $i=1,\dots,n$, choose a plumbing neighborhood $P_i\cong X(p_i,p_i-1)$ of the corresponding $A_{k_i}$-configuration, and set $M_0:=M-\bigcup_{i=1}^n \operatorname{int}P_i$. Then $\partial M_0=\coprod_{i=1}^n -L(p_i,p_i-1)$, where $-L(p_i,p_i-1)$ is the orientation reversal of $L(p_i,p_i-1)$. Also, we have $b_2(M_0)=b_2^+(M_0)=1$. 

Let $\mathfrak{s}_M$ denote the canonical $\operatorname{spin}^c$ structure on $M$, induced by any compatible almost complex structure on $M$. Then $c_1(\mathfrak{s}_M)=-K_M$. For every component sphere $S$ of the $A_{k_i}$-configurations, we have $\langle K_M,[S]\rangle=0$ by the adjunction formula, so we have $K_M|_{P_i}=0\in H^2(P_i;\mathbb{Z})$. Since $K^0_{p_i,p_i-1}=0\in H^2(P_i;\mathbb{Z})$, we must have $\mathfrak{s}_M|_{P_i}=\mathfrak{s}_{K^0_{p_i,p_i-1}} \in \operatorname{Spin}^c(P_i)$. Now consider the following commutative diagram, where the right vertical map is the canonical identification under orientation reversal, and all the other maps are restriction maps.   \[
    \begin{tikzcd}[matrix scale=0.3]
    &&& \operatorname{Spin}^c(M_0) \ar[rrr] &&& \operatorname{Spin}^c(-L(p_i,p_i-1)) \ar[dd, leftrightarrow, "\cong"] \\
        \operatorname{Spin}^c(M) \ar[rrru] \ar[rrrd] \\&&&  \operatorname{Spin}^c(P_i)\ar[rrr] &&& \operatorname{Spin}^c(L(p_i,p_i-1)),
    \end{tikzcd}
    \]
Put $\mathfrak{s}_0:=\mathfrak{s}_M|_{M_0}\in \operatorname{Spin}^c(M_0)$. Applying the diagram to $\mathfrak{s}_M$, and using \eqref{eq:restriction}, we see that $\mathfrak{s}_0|_{-L(p_i,p_i-1)}=\mathfrak{s}_{p_i,p_i-1}(p_i-1)$. 

Since $b_1(M)=0$, $b_2^+(M)=1$, and $b_2^-(M)=N$, we have $K_M^2=2e(M)+3\sigma(M)=9-N$. Hence $c_1(\mathfrak{s}_0)^2=9-N$, since the restriction of $K_M$ to each $P_i$ vanishes. Thus, if we let $\mathfrak{s}_0^-\in \operatorname{Spin}^c(-M_0)$ denote the $\operatorname{spin}^c$ structure corresponding to $\mathfrak{s}_0$ under orientation reversal, then $c_1(\mathfrak{s}_0^-)^2=-(9-N)$. 

Next, let $X:=\coprod_{i=1}^n X(p_i,1)$ and let $\mathfrak{s}_X:=\coprod_{i=1}^n \mathfrak{s}_{K_{p_i,1}}\in \operatorname{Spin}^c(X)$. If $v_i\in H_2(X(p_i,1);\mathbb{Z})$ is the standard generator with dual $v_i^*\in H^2(X(p_i,1);\mathbb{Z})$, then $K_{p_i,1}^2=(p_iv_i^*)^2=-p_i$, and hence $c_1(\mathfrak{s}_X)^2=-\sum_{i=1}^n p_i$. Moreover, $\mathfrak{s}_X|_{L(p_i,1)}=\mathfrak{s}_{K_{p_i,1}}|_{L(p_i,1)}=\mathfrak{s}_{p_i,1}(0)$ by \eqref{eq:restriction}. 

Consider the orientation-reversing diffeomorphism \[
f=\coprod_{i=1}^n f_{p_i,p_i-1}: \partial(-M_0)=\coprod_{i=1}^n L(p_i,p_i-1) \longrightarrow \coprod_{i=1}^n L(p_i,1)=\partial X. 
\]
We have $f^*(\mathfrak{s}_X|_{\partial X})=\mathfrak{s}_0^-|_{\partial(-M_0)}$ by \eqref{eq:pullback}, so there exists a global $\operatorname{spin}^c$ structure $\mathfrak{s}_W$ on  a closed, smooth $4$-manifold $W:=X\cup_f(-M_0)$ such that $\mathfrak{s}_W|_X=\mathfrak{s}_X$ and $\mathfrak{s}_W|_{-M_0}=\mathfrak{s}_0^-$. Note that \begin{align*}
    c_1(\mathfrak{s}_W)^2 & =c_1(\mathfrak{s}_X)^2+c_1(\mathfrak{s}_0^-)^2 =    \sum_{i=1}^n K_{p_i,1}^2 + c_1(\mathfrak{s}_0^-)^2=-\sum_{i=1}^n p_i-(9-N)=-(n+9).
\end{align*}
Since $W$ is a closed, smooth, negative definite $4$-manifold with $b_2(W)=n+1$, Donaldson's diagonalization theorem \cite{Donaldson-1987} implies that $Q_W\cong -\mathbb{Z}^{n+1}$. Let $\mathbf{K}\in -\mathbb{Z}^{n+1}$ be the vector corresponding to (the free part of) $c_1(\mathfrak{s}_W)$ under Poincar\'e duality and the chosen diagonalization. Since the vector $\mathbf{K}$ is characteristic with $\mathbf{K}^2=-(n+9)$, exactly one coefficient has absolute value three and all the others have absolute value one. After a signed permutation of the diagonal basis, we may therefore assume that $\mathbf{K}=3e_1+e_2+\cdots+e_{n+1}$. 

Finally, composing the inclusion $Q_X\hookrightarrow Q_W$ with the chosen diagonalization $Q_W\cong -\mathbb{Z}^{n+1}$ gives a codimension-one lattice embedding \[ 
\iota: Q_X= \bigoplus_{i=1}^n\langle-p_i\rangle \longrightarrow-\mathbb Z^{n+1}.
\]
For each standard generator $v_i$, we have
\[ \mathbf{K}\cdot \iota(v_i)=\langle c_1(\mathfrak{s}_W), v_i\rangle=\langle c_1(\mathfrak{s}_X),v_i\rangle=\langle K_{p_i,1},v_i\rangle=p_i=k_i+1. \]
This proves the proposition.
\end{proof}


\newcommand{\etalchar}[1]{$^{#1}$}

\end{document}